\documentclass[a4paper, 11pt, twoside, dvips]{article}
\usepackage{amsmath}
\usepackage{amsfonts}
\usepackage{float}
\usepackage{amssymb}
\usepackage{amsthm}
\usepackage{amsfonts}
\usepackage{verbatim}
\usepackage{graphicx}
\usepackage{psfrag}
\usepackage[T1]{fontenc}
\usepackage{fancyhdr}
\newtheorem{thm}{Theorem}[section]
\newtheorem{prop}[thm]{Proposition}
\newtheorem{lemma}[thm]{Lemma}
\newtheorem{defn}[thm]{Definition}
\newtheorem{cor}[thm]{Corollary}
\newtheorem{rem}[thm]{Remark}

\newcommand{\df}{\stackrel{\mathrm{def}}{=}}
\newcommand{\R}[1]{\mathbb{R}^{#1}}
\newcommand{\U}[1]{U_{#1,b,L}}
\newcommand{\T}[2]{T^{#1}_{#2}}
\newcommand{\Xt}{X_{t}}
\newcommand{\E}{\mathbb{E}}

\newcommand{\om}{\omega}
\numberwithin{equation}{section}

\begin{document}
\title{Diffusions in Random Environment and Ballistic Behavior} 
\author{Tom Schmitz\\
        Department of Mathematics\\
        ETH Zurich\\
        CH-8092 Zurich\\
        Switzerland\\
        email: schmitz@math.ethz.ch
}
\date{Revised version\\ June 20, 2005}
\maketitle

\bibliographystyle{alpha}
{\bf Abstract:} In this article we investigate the ballistic behavior
of diffusions in random environment. We introduce
conditions in the spirit of $(T)$ and $(T')$ of the discrete setting, cf.
\cite{szn01}, \cite{szn02}, that imply, when $d \geq 2$, a law of large numbers with non-vanishing limiting velocity (which we refer to as 'ballistic behavior') and a central limit theorem with non-degenerate covariance matrix.
As an application of our results, we consider the class of diffusions where the diffusion matrix is the identity, and give a concrete criterion on the drift 
term under which the diffusion in random environment exhibits ballistic 
behavior.
This criterion provides examples of diffusions in random environment with ballistic behavior, beyond what was previously known.\\

{\bf R\'esum\'e:} On \'etudie dans cet article le comportement ballistique de diffusions en milieu al\'eatoire. On montre que certaines conditions $(T)$ et $(T')$, d'abord introduites dans le cadre discret, cf. \cite{szn01}, \cite{szn02}, entra\^inent en dimension sup\'erieure une loi des grands nombres avec une vitesse limite non nulle (ce qu'on appelle 'comportement ballistique'),
et un th\'eor\`eme limite central avec une matrice de covariance non d\'eg\'en\'er\'ee. Pour illustrer ces r\'esultats, on consid\`ere la classe de diffusions 
o\`u la matrice de diffusion est l'identit\'e, et on donne un crit\`ere concret sur la d\'erive qui entra\^ine le comportement ballistique de la diffusion en milieu al\'eatoire. Ce crit\`ere fournit de nouveaux examples de diffusions en milieu al\'eatoire avec comportement ballistique.

\section{Introduction}

The method of ``the environment viewed from the particle'' has played a prominent role in the investigation of random motions in random environment,
see for instance \cite{kip-Var}, \cite{kozlov},  \cite{molchanov}, \cite{olla94}, \cite{olla01}, \cite{papa}, \cite{ras}. In the continuous space-time setting, it applies successfully when one can construct, most often explicitly, 
 an invariant measure for the process of the environment viewed from the particle, 
 which is absolutely continuous with respect to the static measure of 
the random medium, see \cite{deMasi}, \cite{kom}, \cite{kom-krupa02}, \cite{kom-krupa}, \cite{kom-olla-01}, \cite{kom-olla-03}, \cite{landim-olla-yau}, \cite{oel}, \cite{olla94}, \cite{olla01}, \cite{papa}. However, the existence of such invariant measures is hard to prove in the general setting. The case of Brownian motion 
with a random drift which is either incompressible or the gradient of a stationary function, is tractable, see \cite{olla94}, \cite{olla01}. But many examples fall outside this framework, and only recent developments go beyond it, for they require new techniques, see \cite{kom}, 
\cite{kom-krupa02}, \cite{kom-krupa}, \cite{kom-olla-03}.\\ 
Progress has recently been made in the discrete setting for random walks in random environment in higher dimensions, in particular with the help of the renewal-type arguments introduced in Sznitman-Zerner \cite{szn-zer}, see  \cite{bolt-szn-1},  \cite{bolt-szn},  \cite{bolt-szn-zeit}, \cite{com-zeit}, \cite{szn00}, \cite{szn01}, \cite{szn02}, \cite{szn03}, \cite{szn04}, \cite{zeit}. It is natural,
but not straightforward, to try to transpose these results to the continuous space-time setting, and thus propose a new approach to multidimensional diffusions in random environment, when no
invariant measure is a priori known. The first step in this direction was taken up in Shen \cite{shen}, where, in the spirit of Sznitman-Zerner \cite{szn-zer},
certain regeneration times providing a renewal structure are introduced. 
Then a sufficient condition for a 'ballistic' strong law of large numbers 
('ballistic' means that the limiting velocity does not vanish, which we refer to as ballistic behavior) and a central limit theorem governing corrections to the
law of large numbers,
with non-degenerate covariance matrix, is given in terms of these regeneration
times.\\

In this article we show that under condition ($T'$), see (\ref{eq:T'}) for the definition, when $d \geq 2$, the diffusion in random environment
satisfies the aforementioned sufficient condition of Shen \cite{shen}. 
We formulate the rather geometric condition $(T')$ and are able to restate it equivalently in terms of the renewal structure of Shen \cite{shen}, see Theorem
\ref{thm:Tgamma}. With $(T')$ we are then able to derive tail estimates on the first regeneration time which in particular imply the above mentioned 
sufficient condition of Shen \cite{shen}, see Theorem \ref{thm:tail-estimate}.
In the discrete i.i.d. setting, condition ($T'$) was introduced in the work of Sznitman, see \cite{szn01} and \cite{szn02}, and some of our arguments are inspired by
\cite{szn01} and \cite{szn02}.
As an application of our methods, we give concrete examples. In particular,
we recover and extend results of Komorowski and Krupa \cite{kom-krupa}.\\
Before describing our results in more details, let us recall the setting.\\
The {\it random environment} is described by a probability space $(\Omega,\mathcal{A},\mathbb{P})$.  
We assume that there exists a group $\{t_{x}:x \in \R{d}\}$ of transformations on $\Omega$, jointly measurable in $x, \om$, which preserve the probability $\mathbb P$:
\begin{equation}
  \label{eq:stationarity}
  t_x \mathbb P =\mathbb P \,.
\end{equation}
On $(\Omega,\mathcal{A},\mathbb{P})$ we consider bounded measurable functions
$b(\cdot):\Omega \rightarrow \R{d}$ and $\sigma(\cdot):\Omega \rightarrow 
\R{d \times d}$, as well as two finite constants $\bar{b},~\bar{\sigma}>0$ 
such that for all $\om \in \Omega$
\begin{eqnarray}
\label{eq:b-sigma-bound}
\left| b(\om) \right| \leq \bar{b}, \quad \left| \sigma(\om) \right| \leq 
\bar{\sigma},
\end{eqnarray}
where $|\cdot|$ denotes the Euclidean norm for vectors resp. for square 
matrices.
We write  
\begin{eqnarray*}  
b(x,\om)=b(t_{x}(\om)), \quad \sigma(x,\om)=\sigma(t_{x}(\om)).
\end{eqnarray*} 
We further assume that $b(\cdot,\om)$ and $\sigma(\cdot,\om)$ are
Lipschitz continuous, i.e. there is a constant $K>0$ such that for
all $\om \in \Omega,~x,y \in \R{d}$,
\begin{eqnarray}
\label{eq:Lipschitz}
|b(x,\om)-b(y,\om)|+|\sigma(x,\om)-\sigma(y,\om)| \leq K|x-y|.
\end{eqnarray}
$\sigma \sigma^{t}(x,\om)$ is uniformly elliptic, i.e. there is a constant
$\nu > 0$ such that for all $\om \in \Omega,~x,y \in \R{d}$,
\begin{eqnarray}
\label{eq:elliptic}
\frac{1}{\nu}|y|^{2} \leq |\sigma^{t}(x,\om)y|^{2} \leq \nu |y|^{2},
\end{eqnarray}
where $\sigma^{t}$ denotes the transposed matrix of $\sigma$.
For a Borel subset $F \subset \R{d}$, we define the $\sigma$-field
generated by $b(x,\om),~\sigma(x,\om)$, for $x \in F$ by
\begin{equation}
\label{eq:sigma-field}
\mathcal{H}_{F} \df \sigma\{b(x,\cdot), \sigma(x,\cdot):x \in F\},
\end{equation}
and assume finite range dependence: there is an $R > 0$ such that for all Borel subsets $F, F' \subset \R{d}$ with $d(F,F') \df \inf\{|x-x'|: x \in F, x' \in F' \}>R$,
\begin{equation}
\label{eq:R-separation}
\mathcal{H}_{F} \text{ and }  \mathcal{H}_{F'} \text{ are $\mathbb P$-independent}. 
\end{equation}
We denote by $(C(\R{}_{+},\R{d}),\mathcal{F},W)$ the canonical Wiener space,
and with $(B_{t})_{t \geq 0}$ the $d$-dimensional Brownian motion
(which is independent from $(\Omega,\mathcal{A},\mathbb{P})$).
The diffusion process in the random environment $\om$ is described by the 
family of laws $(P_{x,\om})_{x \in \R{d}}$ (we call them the \emph{quenched} laws)
on $(C(\R{}_{+},\R{d}),\mathcal{F})$ of the solution of the stochastic
differential equation 
\begin{eqnarray}
\label{eq:SDE}
 \begin{cases}
  dX_t=\sigma(X_t,\om)dB_t+b(X_t,\om)dt,\\
  X_{0}=x, \quad  x \in \R{d}, ~\om \in \Omega.
 \end{cases}
\end{eqnarray}
The second order linear differential operator associated to the stochastic
differential equation (\ref{eq:SDE}) is given by:
\begin{equation}
  \label{eq:diff-operator}
  \mathcal{L_\om} \df \frac{1}{2}\sum_{i,j=1}^d a_{ij}(x,\om)\frac{\partial^2}{\partial  x_i\partial x_j  }+\sum_{j=1}^d b_j(x, \om)\frac{\partial}{\partial x_j}\,.
\end{equation}
To restore some stationarity to the problem, it is convenient to introduce the \emph{annealed} laws $P_{x}$, which are defined as the 
semi-direct products:
\begin{equation}
\label{eq:annealed}
P_{x} \df \mathbb{P} \times P_{x,\om},~~ \mathrm{for}~ x \in \R{d}.
\end{equation}
Of course the Markov property is typically lost under the annealed laws. \\

Let us now explain the purpose of this work.
The main object is to introduce sufficient conditions for ballistic behavior of the diffusion in random environment when $d \geq 2$.
These conditions are expressed in terms of another condition $(T)_{\gamma}$
which is defined as follows.
Consider, for $|l|=1$ a unit vector of $\R{d}$, $b, L>0$, the slabs 
\begin{equation*}
U_{l,b,L} \df \{x \in \R{d}:-bL < x \cdot l <L\}.
\end{equation*}
We say that \emph{condition $(T)_{\gamma}$} holds relative to
$l \in S^{d-1}$, in shorthand notation $(T)_\gamma |\,l$, if for all $l' \in S^{d-1}$ in a neighborhood of $l$, and for all $b>0$,
\begin{equation}
\label{eq:Tgamma}
\limsup_{L \to \infty} L^{-\gamma} \log{P_{0}[X_{\T{}{\U{l'}}}\cdot l' < 0]}<0
,\end{equation}
where $\T{}{\U{l}}$ denotes the exit time of $X_\cdot$ out of the slab $U_{l,b,L}$, see (\ref{eq:exit-time}) for the definition.\\
The aforementioned sufficient conditions for ballistic behavior are then condition $(T)$ relative to the direction $l$, in shorthand notation $(T)|l$, which refers to the case where
\begin{equation}
\label{eq:T}
(\ref{eq:Tgamma}) \textrm{ holds for } \gamma = 1\,,
\end{equation}
or the weaker condition $(T')$ relative to the direction $l$, in shorthand notation $(T')|l$, which refers to the case where 
\begin{equation}
\label{eq:T'}
(\ref{eq:Tgamma}) \textrm{ holds for all } \gamma \in (0,1)\,.
\end{equation} 
Clearly $(T)$ implies $(T')$ which itself implies
$(T)_{\gamma}$ for all $\gamma \in (0,1)$. We expect these conditions all
to be equivalent, cf. Sznitman \cite{szn02}, \cite{szn04}, however this remains an open question.
The conditions $(T)$ and $(T')$ are not effective conditions which can be checked by direct inspection of the environment restricted to a bounded
domain of $\R{d}$. 
In the discrete i.i.d. setting, Sznitman \cite{szn02} 
proved the equivalence between a certain effective criterion and condition $(T')$. 
With the help of the effective criterion he also proved that $(T)_\gamma$ and
$(T')$ are equivalent for $\frac{1}{2}< \gamma <1$.
We believe that a similar effective criterion holds in the continuous setting, and it is in the spirit of this 
belief that we formulate all our results in Section \ref{sec:condition-T} and \ref{sec:tail-estimate} in terms of
condition $(T')$ resp. $(T)_\gamma$.  
Later, in Section \ref{sec:examples}, we verify the stronger condition ($T$) 
for a large class of examples.\\
In Theorem \ref{thm:Tgamma} we show that the definition of condition $(T)_\gamma|l$, see (\ref{eq:Tgamma}), which is of a rather geometric nature, has an equivalent formulation in terms of transience of the diffusion in direction $l$ and a stretched exponential control of 
the size of the trajectory up to the first regeneration time $\tau_{1}$ (see subsection \ref{subsec:regeneration} for the precise definition):
\begin{align}
& P_{0}-a.s. \lim_{t \to \infty} X_t \cdot l = \infty\,, 
\label{eq:transience-0}
\\
&\text{and for some } \mu > 0,~~ \hat{E_{0}}[\,\exp\{\mu \sup_{0 \leq t \leq \tau_{1}}|\Xt|^{\gamma}\}] < \infty\,.
\label{eq:integrability}
\end{align}
Following Shen \cite{shen}, the successive regeneration times $\tau_{k}, k \geq 1$, are defined on an enlarged
probability space which is obtained by adding some suitable auxiliary i.i.d.~Bernoulli variables, cf. subsection \ref{subsec:regeneration}. The quenched measure on the enlarged space, which couples the trajectories to the Bernoulli variables,
is denoted by $\hat{P}_{x,\om}$, and $\hat{P}_{x}$ refers to the annealed
measure $\mathbb{P}\times \hat P_{x,\om}$, cf. subsection \ref{subsec:coupling}. 
Loosely speaking, the
first regeneration time $\tau_{1}$ is the first integer time where the 
diffusion process in random environment reaches a local maximum in a given
direction $l \in S^{d-1}$, some auxiliary Bernoulli variable takes value one,
and from then on the diffusion process never backtracks. \\
The strategy of the proof of the above mentioned equivalence statement is similar to that of the analogue statement in the discrete i.i.d. setting, see Sznitman \cite{szn02}. Nevertheless, changes
appear in several places, due among others to the fact that the 
regeneration time $\tau_{1}$ is more complicated than in the discrete setting.\\
 
Theorem \ref{thm:Tgamma} is very useful because conditions (\ref{eq:Tgamma}) and (\ref{eq:integrability}) have different flavours.
Condition (\ref{eq:integrability}) is especially useful when studying asymptotic
properties of the diffusion process, whereas (\ref{eq:Tgamma}) is more
adequate to construct examples.\\

Together with the crucial renewal
property (see Theorem \ref{thm:renewal}) induced by the regeneration times $\tau_{k}$, $k \geq 1$, the formulation (\ref{eq:integrability}) is
instrumental in showing that under $(T')$, and when $d \geq 2$,
\begin{equation}
\label{eq:tail}
\limsup_{u \rightarrow \infty}~ (\log u)^{-\alpha} \log \hat{P}_{0}[\tau_{1} > u]<0, \quad \mathrm{for}
\quad \alpha < 1+\frac{d-1}{d+1}, 
\end{equation}
see Theorem \ref{thm:tail-estimate}.
The proof again uses a strategy close to the proof in the discrete case, see 
Sznitman \cite{szn01}. We prove a seed estimate, see Lemma \ref{lemma:seed-estimate}, which is then propagated to the right scale by performing
a renormalisation step, see Lemma \ref{lemma:renormalisation}. Interestingly enough, 
we do not require condition $(T')$ to prove the renormalisation lemma. 
\\
Under the assumption of (\ref{eq:transience-0}) and the finiteness of the first and the second moment of $\tau_1$, the Theorems 3.2 and 3.3 in Shen \cite{shen} imply that:
\begin{equation}
\label{eq:lln}
P_{0}-\text{a.s.}, \quad \frac{X_{t}}{t} \rightarrow v, \quad v
\neq 0,  \text{ deterministic, with } v \cdot l >0\,,  
\end{equation}
\begin{equation}
\label{eq:clt}
\begin{aligned}
&\text{and under $P_{0}$, $B_{\cdot}^{s}=\tfrac{X_{s \cdot}-s\cdot v}{\sqrt{s}}$
 converges in law on $C(\R{}_{+},\R{d})$, as $s \to \infty$, to a}\\
&\text{Brownian motion $B_{\cdot}$ with non-degenerate deterministic covariance matrix}.
\end{aligned}
\end{equation}
Hence, when condition $(T')$ holds, and $d \geq 2$, Theorem  \ref{thm:tail-estimate}, see also (\ref{eq:tail}), yields a ballistic law of large numbers and a central limit theorem governing corrections to the law of large numbers.
Incidentally let us mention that as in the discrete setting, cf. Sznitman \cite{szn02}, \cite{szn04}, condition $(T')$ is a natural
contender for the characterisation of ballistic diffusions in random environment when $d \geq 2$. However at present there are no rigorous results in that direction.\\

As an application of our methods, we provide a rich class of examples
exhibiting ballistic behavior.
We first consider the case where, for some $l \in S^{d-1}$ and all $\om \in \Omega$, all $x \in \R{d}$, $b(x, \om)\cdot l $ remains uniformly positive, and show in Proposition \ref{prop:non-nestling} that condition $(T)|l$ holds.
Hence we recover and extend the main result of Komorowski and Krupa \cite{kom-krupa} (which only asserts (\ref{eq:lln}) when $\sigma = Id$).\\
Then we consider the case where $\sigma$ in (\ref{eq:SDE}) is the identity. We prove in Theorem \ref{thm:examples} that, when $d \geq 1$, there is a constant $c_e(\bar b, K, d, R)>0$ such that, for $l \in S^{d-1}$,
\begin{equation}
\label{eq:examples}
\E[(b(0,\omega) \cdot l)_{+}]>c_e~ \E[(b(0,\omega) \cdot l)_{-}]
\end{equation}
implies condition $(T)|l$ (and hence condition $(T')|l$).
Clearly, when $\sigma = Id$, the result of Proposition \ref{prop:non-nestling} is included in Theorem \ref{thm:examples}.  
Note that Theorem \ref{thm:examples} covers additional situations where
$b(0,\om) \cdot l$ changes sign in every unit direction $l$. This provides
 new examples of ballistic diffusions in random environment.
More details are included in remark \ref{rem:examples} at the end of Section \ref{sec:examples}.\\
To prove Theorem \ref{thm:examples}, we verify the geometric formulation 
(\ref{eq:T}) of condition ($T$).
However it is a difficult task to compute the exit 
distribution of the diffusion out of large slabs under $P_{0}$, since 
the Markov property is lost under $P_{0}$.
In the spirit of Kalikow \cite{kalikow}, we restore a Markovian character
to the exit problem by virtue of Proposition \ref{prop:exit}. With the help of
Proposition \ref{prop:exit}, we show that condition $(T)$ is implied by a certain condition ($K$), see (\ref{eq:K}), which has a similar flavor as Kalikow's condition in the discrete i.i.d. setting, see Sznitman and Zerner \cite{szn-zer}.
 The proof of Theorem \ref{thm:examples} is then carried out by checking condition ($K$).
These steps are similar in spirit to the strategy used in the discrete setting, cf. lecture 5 of \cite{bolt-szn}.
However, difficulties arise in the continuous space-time framework.\\ 
Let us now describe the organisation of this article.\\
In Section \ref{sec:recall}, we recall the coupling construction which leads to the measures $\hat P_{x,\om}$ resp. $\hat P_x$, cf. Proposition \ref{prop:coupling}. On this new probability space one constructs the regeneration times $\tau_{k}, ~k \geq 1$, which provide the crucial renewal structure, cf. Theorem \ref{thm:renewal}. These results have been obtained in Shen \cite{shen}; we recall them for the convenience of the reader.\\
In Section \ref{sec:condition-T}, we prove the equivalence of (\ref{eq:Tgamma}) and (\ref{eq:transience-0}), (\ref{eq:integrability}), see Theorem \ref{thm:Tgamma}. \\
In Section \ref{sec:tail-estimate}, we show (\ref{eq:tail}) under the assumption of condition $(T')$, see Theorem \ref{thm:tail-estimate}.
Proposition \ref{prop:trap} highlights the importance of large deviation
controls of the exit probability of large slabs.
The renormalisation step is carried out in Lemma \ref{lemma:renormalisation}, and a seed estimate is provided in Lemma \ref{lemma:seed-estimate}.\\
In Section \ref{sec:examples}, we show that condition $(T)$ (in the geometric formulation (\ref{eq:Tgamma}))
holds either under the assumption of the uniform positivity of $b(x,\om)\cdot l$ for some unit vector $l$ and all $\om \in \Omega$, all $x \in \R{d}$, or
under the assumption of $\sigma=Id$ and (\ref{eq:examples}). \\
In the Appendix, we provide some results on continuous local martingales and Green functions, that we use throughout this article.\\

{\bf Convention on constants}
Unless otherwise stated, constants only depend on the quantities $\nu, \bar{b},\bar{\sigma},K, R, d, \gamma$. In particular they are independent of the environment $\om$. Generic positive constants are denoted by $c$. 
Dependence on additional parameters appears in the notation. For example, $c(p,L)$ means that the constant $c$ depends on $p$ and $L$ {\it and on} $\nu, \bar{b},\bar{\sigma},K, R, d, \gamma$.
When constants or positive numbers are not numerated, their value may change from line to line.\\ 

{\bf Acknowledgement:} Let me thank my advisor Prof. A.-S. Sznitman
for introducing me to the subject and for his advice during the completion of this work. I also want to thank Lian Shen for his help and numerous discussions.  
\section{The Regeneration Times and the Renewal Structure}
\label{sec:recall}
In this section, we recall the definition of the coupled measures $\hat{P}_{x,\om}$ (resp. $\hat{P}_{x}$) and of the regeneration times $\tau_{k}$, $k \geq 1$, given in Shen \cite{shen}. We then cite the resulting renewal structure, see Theorem \ref{thm:renewal}.
For the proofs or further details, we refer the reader to Shen \cite{shen}.
\subsection{Notation}
We introduce some additional notation. For $x \in \R{d}$, $d \ge 1$, we let $B_{r}(x)$ denote the open 
Euclidean ball with radius $r$ centered in $x$. For $U \subseteq \R{d}$, 
we denote with $\bar U$ its closure, with diam$(U) \df \sup\{|x-y|:x,y \in U\}$ its diameter, and, for measurable $U$, with $|U|$ its Lebesgue measure. A domain stands for a connected open subset of $\R{d}$.
For $x \in \R{}$, we define $\lfloor x \rfloor \df \sup\{k \in \mathbb Z:k \leq x\}$ and $\lceil x \rceil \df \inf\{k \in \mathbb Z:k \geq x\}$. For a discrete set $A$, we denote with $\#A$ its cardinality.
For an open set $U$ in $\R{d}$ and $u \in \R{}$ we define the 
$(\mathcal{F}_t)_{t \geq 0}$-stopping times (($\mathcal{F}_t)_{t \geq 0}$ 
denotes the canonical right-continuous filtration on $(C(\mathbb{R}_{+}, \R{d}), \mathcal{F}))$:\\
the exit time from $U$,
\begin{equation}
\label{eq:exit-time}
T_U \df \inf{\{t \geq 0:\Xt \notin U\}},
\end{equation}
and the entrance times into the half-spaces $\{x \cdot l \geq u\}$ resp. 
 $\{x \cdot l \leq u\}$,
\begin{align}
\label{eq:stopping-times}
\begin{split}
&\T{l}{u} \df \inf{\{t \geq 0:X_t \cdot l \geq u\}}, \\
&\widetilde{T}^{l}_{u} \df \inf{\{t \geq 0:X_t \cdot l \leq u\}}. 
\end{split}
\end{align}
We define as well the maximal value of the process $(X_{s}\cdot l)_{s \geq 0}$ till time t,
\begin{equation}
\label{eq:maximum}
M(t) \df \sup{\{X_{s} \cdot l: 0 \leq s \leq t \}},
\end{equation}
and the first return time of the process $(X_s \cdot l)_{s \geq 0}$
to the level $-R$ relative to the starting point, as well as its rounded value,
\begin{equation}
  \label{eq:return}
  J \df \inf\{t \geq 0: (X_t-X_0)\cdot l \leq -R\}\, , \quad D \df \lceil J \rceil \, .
\end{equation}
\subsection{The coupled measures}
\label{subsec:coupling}
We need further notations. We let $l$ be a fixed unit vector, and 
\begin{equation}
  \label{eq:2-subsets}
  U^x\df B_{6R}(x+5Rl)\,,\quad B^x\df B_R(x+9Rl)\,.
\end{equation}
We denote by $\lambda_j$ the canonical coordinates on $\{0,1\}^\mathbb N$. 
Further, let $\mathcal S_m\df \sigma\{\lambda_0,\cdots,\lambda_m\}$\,, $m\in\mathbb N$, denote the canonical filtration on $\{0,1\}^\mathbb N$ generated by $(\lambda_m)_{m\in\mathbb N}$ and $\mathcal S \df \sigma\big\{\bigcup_m \mathcal S_m\big\}$ be the canonical $\sigma$-algebra. We also write for $t\geq 0$:
\begin{equation}
  \label{eq:filtrationZ}
  \mathcal Z_t\df \mathcal F_t\otimes\mathcal S_{\lceil t\rceil}\,,\quad \mathcal Z \df \mathcal F\otimes\mathcal S=\sigma\Big\{\bigcup_{m\in\mathbb N}\mathcal Z_m\Big\}\,.
\end{equation}
We also consider the shift operators $\big\{\theta_m:m\in\mathbb N\big\}$, with $\theta_m: 
\left(C(\mathbb R_+,\R{d})\times \{0,1\}^\mathbb N, \mathcal Z\right)\to \left(C(\R{}_+,\R{d})\times \{0,1\}^\mathbb N,\mathcal Z\right)$, such that
\begin{equation}
  \label{eq:shift}
  \theta_m\left(X_\cdot,\lambda_\cdot\right)=\left(X_{m+\cdot},\lambda_{m+\cdot}\right)\,. 
\end{equation}
Then from Theorem 2.1 in Shen \cite{shen}, one has the following measures, coupling the diffusion in random environment with a sequence of Bernoulli variables:
\begin{prop}
\label{prop:coupling}
There exists $p>0$, such that for every $\om\in\Omega$ and $x\in\R{d}$,  there exists a probability measure $\hat P_{x,\om}$ on $\big(C(\R{}_+,\R{d})\times\{0,1\}^\mathbb N,\mathcal Z\big)$ depending measurably on $\om$ and $x$, such that
  \begin{enumerate}
  \item Under $\hat P_{x,\om}$, $(X_t)_{t\geq 0}$ is $P_{x,\om}$-distributed, and the $\lambda_m$, $m\geq 0$, are i.i.d. Bernoulli variables with success probability $p$.
  \item  For $m \geq 1$, $\lambda_m$ is independent of 
$\mathcal F_m\otimes\mathcal S_{m-1}$ under $\hat P_{x,\om}$. Conditioned on $\mathcal Z_m$, $X_\cdot\circ\theta_m$ has the same law as $X_\cdot$ under $\hat P^{\lambda_m}_{X_m,\om}$, where for $y \in \R{d}$, $\lambda \in \{0,1\}$, $\hat P^{\lambda}_{y,\om}$ denotes the law 
$\hat P_{y,\om}[\;\cdot\;|\lambda_0=\lambda]$.
  \item $\hat P^{1}_{x,\om}$ almost surely, $X_{s}\in U^{x}$ for $s\in[0,1]$ (recall (\ref{eq:2-subsets})).
  \item Under $\hat P^{1}_{x,\om}$, $X_{1}$ is uniformly distributed on $B^x$ (recall (\ref{eq:2-subsets})).
  \end{enumerate}
\end{prop}
We then introduce the new annealed measures on $\big(\Omega\times C(\R{}_+,\R{d})\times\{0,1\}^\mathbb N,\mathcal A\otimes\mathcal Z\big)$:
\begin{equation}
  \label{eq:coupled-measures}
  \hat P_{x} \df \mathbb P\times \hat P_{x,\om} \quad\text{and}\quad \hat E_x
  \df \mathbb E\times \hat E_{x,\om}\,.
\end{equation}

\subsection{The Regeneration Times $\tau_k$ and the Renewal Structure}
\label{subsec:regeneration}
To define the first regeneration time $\tau_1$, we introduce a sequence of integer-valued $(\mathcal Z_t)_{t\geq 0}$-stopping times $N_k$, $k \ge 1$, such that, at these times, the  Bernoulli variable takes the value one, and the process $(X_t \cdot l)_{t\geq 0}$ in essence reaches a new maximum.
Proposition \ref{prop:coupling} now shows that {\it for every environment $\om \in \Omega$}, the position of the diffusion at time $N_{k+1}$ is uniformly distributed on the ball $B^{X_{N_k}}$ under $\hat P_{0,\om}$. We define $\tau_1$ as the first $N_k+1$ such that, after time $N_k+1$, 
the process $(X_t \cdot l)_{t \ge 0}$ never goes below the level $X_{N_k+1} \cdot l -R$. 
In essence, the distance between the positions $X_{\tau_1-1}$ and $X_{\tau_1}$ 
is large enough to obtain, in view of finite range dependence, independence
of the parts
of the trajectory $(X_t-X_0)_{t \le \tau_1-1}$ and $(X_{\tau_1+t}-X_{\tau_1})_{t \ge 0}$ under $\hat P_0$, so that the diffusion regenerates at time $\tau_1$
under $\hat P_0$. 
We define the regeneration times $\tau_k$, $k \ge 2$, in an iterative fashion, and we provide the renewal structure  in Theorem \ref{thm:renewal}.\\
In fact, the precise definition of $\tau_{1}$ relies on several sequences of stopping times.
First, for $a>0$, introduce the $(\mathcal F_t)_{t\geq 0}$-stopping times $V_k(a)$, $k\geq 0$, (recall $M(t)$ in (\ref{eq:maximum}) and $T_u$ in (\ref{eq:stopping-times})):
\begin{equation}
  \label{eq:V(a)}
  V_0(a)\df T_{M(0)+a}\,, \, V_{k+1}(a)\df T_{M(\lceil V_k(a)\rceil)+R}\,.
\end{equation}
In view of the Markov property, see point 2. of Proposition \ref{prop:coupling}, we want the stopping times $N_k(a), k\geq 1$, to be integer-valued. Therefore we introduce in an intermediate step the (integer-valued) stopping times $\tilde N_k(a)$ where the process $X_t\cdot l$ essentially reaches a maximum:
\begin{equation}
\begin{cases}
  \displaystyle \tilde N_1(a)\df \inf\Big\{\lceil V_k(a)\rceil : k\geq 0, \sup_{s\in[V_k,\lceil V_k\rceil]}|l\cdot(X_s-X_{V_k})|< \tfrac{R}{2}\Big\}\,,\\
        \tilde N_{k+1}(a)\df \tilde N_1(3R)\circ\theta_{\tilde N_k(a)}+\tilde N_k(a)\,,\ k\geq 1\,,
\end{cases}
\end{equation}
(by convention we set $\tilde N_{k+1}=\infty$ if $\tilde N_k=\infty$).
In the spirit of the comment at the beginning of this subsection, we define the $(\mathcal Z_t)_{t\geq 0}$-stopping time $N_1$ as
\begin{equation}
 \label{eq:N(a)}
   N_1(a)\df \inf\left\{\tilde N_k(a): k\geq 1, \lambda_{\tilde N_k(a)}=1\right\},\,\,\,N_1 \df N_1(3R),
\end{equation}
as well as the $(\mathcal Z_t)_{t\geq 0}$-stopping times
\begin{equation}
  \label{eq:time-S}
  \begin{cases}
    S_1\df N_1+1\,,\\    
    R_1\df S_1+D\circ\theta_{S_1}\,.
  \end{cases}
\end{equation}
The $(\mathcal Z_t)_{t\geq 0}$-stopping times $N_{k+1}$, $S_{k+1}$ and $R_{k+1}$ are defined in an iterative fashion for $k \geq 1$: 
\begin{equation}
  \label{eq:N-S-R}
  \begin{cases}
  N_{k+1}\df R_{k}+N_1(a_k)\circ\theta_{R_k} \text{ with } 
  a_k\df M(R_k)- X_{R_k} \cdot l+R\ge R\,,\\
  S_{k+1}\df N_{k+1}+1\,,\\
    R_{k+1}\df S_{k+1}+D\circ\theta_{S_{k+1}}\,                              
  \end{cases}
\end{equation}
(the shift $\theta_{R_k}$ is {\em not} applied to $a_k$ in the above definition).\\
Notice that for all $k \geq 1$, the $(\mathcal Z_t)_{t\geq 0}$-stopping times $N_k$, $S_k$ and $R_k$ are integer-valued, possibly equal to infinity, and we have $1\leq N_1\leq S_1\leq R_1\leq N_2\leq S_2\leq R_2\cdots\leq\infty$.\\
The first {\em regeneration time} $\tau_1$\, is defined, as in \cite{szn-zer}, by
\begin{equation}
  \label{eq:regeneration-time}
  \tau_1\df \inf\{S_k : S_k<\infty,\, R_k=\infty\}\leq \infty\;.
\end{equation}
We define the sequence of random variables $\tau_k$, $k \geq 1$, iteratively on the event  $\{\tau_1<\infty\}$, by viewing $\tau_k$ as a function of $(X_{\cdot}, \lambda_{\cdot})$:
\begin{equation}
  \label{eq:tau_k}
    \tau_{k+1}\big(( X_\cdot,\lambda_\cdot)\big)
  \df \tau_1\big(( X_\cdot,\lambda_\cdot)\big)+\tau_k\big(( X_{\tau_1+\cdot}- X_{\tau_1}, \lambda_{\tau_1+\cdot})\big), \ k\geq 1,
\end{equation}
and set by convention $\tau_{k+1}=\infty$ on $\{\tau_k=\infty\}$. Observe that for each $k \geq 1$, $\tau_k$ is either infinite or a positive integer.
 By convention, we set $\tau_0=0$. 
The random variables $\tau_k$, $k \geq 0$, provide a renewal structure, see also Theorem 2.5 in Shen \cite{shen}, which will be crucial in the proof of Theorem \ref{thm:Tgamma}.
\begin{thm}[Renewal Structure]
\label{thm:renewal}
Assume that $\hat P_0{}$-a.s., $\tau_1<\infty$. Then under the measure $\hat P_0$, the random variables $Z_k\df \left(X_{(\tau_k+\cdot)\wedge(\tau_{k+1}-1)}-X_{\tau_k};\, X_{\tau_{k+1}}-X_{\tau_k};\, \tau_{k+1}-\tau_k\right)$, $k\geq 0$, are independent. Furthermore, $Z_k$, 
$k\geq 1$, under $\hat P_0$, have the distribution of $Z_0=\left(X_{\cdot\wedge(\tau_{1}-1)}-X_0;\, X_{\tau_{1}}-X_0;\, \tau_{1}\right)$ under $\hat P_0[\;\cdot\;|D=\infty]$. 
\end{thm}
The following Proposition is also established in \cite{shen} (see Lemma 2.3 and Proposition 2.7 therein):
\begin{prop}
  \label{prop:equiv}
  $\hat P_0$-a.s. $\tau_1<\infty$ if and only if $P_0$-a.s. $\lim_{t\to
  \infty} X_t\cdot l=\infty$. Furthermore
  $\hat P_0$-a.s. $\tau_1<\infty$ implies $P_0[D=\infty]>0$ (recall the definition of $D$ in (\ref{eq:return})).
\end{prop}

\section{Equivalent Formulations of Condition $(T)_\gamma$}
\label{sec:condition-T}
In this section, we provide an equivalent formulation of the condition $(T)_\gamma |l$, cf. (\ref{eq:Tgamma}), in terms of a stretched exponential estimate on the size of the trajectory $\Xt,0 \leq t \leq \tau_{1}$.
\begin{thm}
\label{thm:Tgamma}
Let $l \in S^{d-1}, 0< \gamma \leq 1$. One has the equivalence
\begin{align}
&\bullet (T)_{\gamma}|l 
\label{eq:geometric-Tgamma} \\
&\bullet P_{0}-a.s. \lim_{t \to \infty} \Xt \cdot l = \infty\,,\text{and
for some }\mu > 0,\, \hat{E_{0}}[\exp\{{\mu \sup_{0 \leq t \leq \tau_{1}}|\Xt|^{\gamma}\}}] < \infty\,.
\label{eq:integrability-2} 
\end{align}
\end{thm} 

\subsection{The Proof of (\ref{eq:geometric-Tgamma}) $\Rightarrow$ 
(\ref{eq:integrability-2})} 
Let us first show that 
\begin{equation}
\label{eq:transience-1}
P_{0}-a.s.\, \lim_{t \to \infty} \Xt \cdot l = \infty\,.
\end{equation} 
We choose an orthonormal basis 
$(f_{i})_{1 \leq i \leq d}$ of $\R{d}$ with $f_{1}=l$. By definition of condition $(T)_\gamma|l$, there are unit vectors $l_{i,+},\, l_{i,-}$
 in $\R{}f_{1}+ \R{}f_{i}$, $2 \leq i \leq d$, such that:
\begin{equation*}
l_{i,\pm} \cdot f_{1} > 0, \quad l_{i,+} \cdot f_{i} > 0,
\quad l_{i,-} \cdot f_{i} < 0, 
\end{equation*}
and, for $l'=l,\, l_{i,+},\, l_{i,-},\, 2 \leq i \leq d,\,b>0$,
\begin{equation}
\label{eq:Tgamma-l'}
 \limsup_{L \to \infty} L^{-\gamma} \log{P_{0}[X_{\T{}{\U{l'}}}\cdot l' < 0]}<0\,. 
\end{equation}
Consider the open set 
$\mathcal{D} \df \{x \in \R{d},\, |x \cdot l|<1,\, x \cdot l_{i,\pm}>-1, \,  2 \leq i \leq d \}.$
$\mathcal{D}$ is a bounded set, hence we can find numbers $a_{i,\pm}>0 \, ,
\, 2 \leq i \leq d$, such that
\[
\mathcal{D} \subseteq \{x \in \R{d}:\, x \cdot l_{i,\pm}<
a_{i,\pm} \, , 2 \leq i \leq d \}\,.
\]
Since $(T)_{\gamma}$ holds relative to $l$ and  $l_{i,\pm}$, 
$2 \leq i \leq d $, writing 
\begin{eqnarray*}
    P_{0} [T_{L\mathcal{D}} < T^{l}_{L}]\,
\leq \,P_{0} [\tilde{T}^{l}_{-L} < T^{l}_{L}] ~+~ \sum_{i=2}^{d}
       P_{0} [\tilde{T}^{l_{i,+}}_{-L} < T^{l_{i,+}}_{La_{i,+}}] ~+~
       \sum_{i=2}^{d}
       P_{0} [\tilde{T}^{l_{i,-}}_{-L} < T^{l_{i,-}}_{La_{i,-}}]\,,
\end{eqnarray*}
we find by (\ref{eq:Tgamma-l'}) that
\begin{equation}
\label{eq:exit-D}
\limsup_{L \to \infty} L^{-\gamma} \log{P_{0}[T_{L\mathcal{D}} < T^{l}_{L}]} <0\,.
\end{equation}
Since $P_{0} [T^{l}_{L} = \infty] \leq P_{0} [T_{L\mathcal{D}} < T^{l}_{L}]$, and the left-hand side increases with $L$, (\ref{eq:exit-D}) implies 
that $P_{0}-a.s. \, \limsup_{t \to \infty} \Xt \cdot l = \infty\,$.
As a next step we observe that 
\begin{equation}
\label{eq:transience}
\limsup_{L \to \infty} L^{-\gamma} \log{P_{0}[\tilde{T}_{\frac{L}{2}}^{l}
\circ \theta _{T_{L}^{l}} < T_{\frac{4L}{3}}^{l}\circ \theta _{T_{L}^{l}}]} <0\,.
\end{equation}
Indeed:
\begin{equation}
\label{eq:transience-3}
     P_{0}[\tilde{T}_{\frac{L}{2}}^{l}\circ \theta _{T_{L}^{l}} < 
       T_{\frac{4L}{3}}^{l}\circ \theta _{T_{L}^{l}}]
\leq  P_{0} [\T{}{L\mathcal{D}} < \T{l}{L}] + 
        P_{0}[\tilde{T}_{\frac{L}{2}}^{l}\circ \theta _{T_{L}^{l}} < 
        T_{\frac{4L}{3}}^{l}\circ \theta _{T_{L}^{l}},
        \T{}{L\mathcal{D}}= \T{l}{L}]\,,  
\end{equation}
and by (\ref{eq:exit-D}) we only need to estimate the second term on the
right-hand side of (\ref{eq:transience-3}). We define
\begin{equation*}
\partial_{+}\mathcal{D} \df \{x \in \partial \mathcal{D}:~x \cdot l =1\},
\end{equation*}
and let $(B_{1}(x_{i}))_{i \in I}$, $x_{i} \in \partial_{+}L\mathcal{D}$, $I$
a finite set with cardinality growing polynomially in $L$, be a cover of
$\partial_{+}L\mathcal{D}$ by unit balls, see above (\ref{eq:stopping-times}) for the notation.
It follows from the strong Markov property and the stationarity of the measure
$\mathbb{P}$ that
\begin{equation}
\label{eq:P}
\begin{aligned}
&P_{0}[\tilde{T}_{\frac{L}{2}}^{l}\circ \theta _{T_{L}^{l}} < 
      T_{\frac{4L}{3}}^{l}\circ \theta_{T_{L}^{l}},
      T_{L\mathcal{D}}= T^{l}_{L}]  
\leq \sum_{i \in I}\mathbb{E}\Big[
       E_{0,\omega}[P_{X_{\T{l}{L}},\omega}[\tilde{T}_{\frac{L}{2}}^{l} <
                  T_{\frac{4L}{3}}^{l}],
                  X_{\T{l}{L}} \in B_1(x_i)]\Big]\\ 
\leq& \sum_{i \in I}\mathbb{E}\left[ \sup{_{x \in B_1(x_i)}P_{x,\omega}
        [\tilde{T}_{\frac{L}{2}}^{l} <T_{\frac{4L}{3}}^{l}]}\right]
= \sum_{i \in I}\mathbb{E}\left[\sup{_{x \in B_{1}(0)}P_{x,\omega}
        [\tilde{T}_{-\frac{L}{2}}^{l} <T_{\frac{L}{3}}^{l}]}\right].
\end{aligned}
\end{equation}
For large enough $L$, it follows from the strong Markov property that for all $\om \in \Omega$, 
\begin{equation}
  \label{eq:strong-Markov}
  \text{the function } x \mapsto  P_{x, \om}[\tilde{T}_{-\frac{L}{2}}^{l} < T_{\frac{L}{3}}^{l}] \text{ is $\mathcal L_\om$-harmonic on $B_3(0)$},
\end{equation}
see for instance \cite{kar-shr} p.364f.
Harnack's inequality (see \cite{gil-tru} p.250) states that there is a constant $c_H>1$ such that for all $\mathcal L_\om$-harmonic functions $u$ on $B_{3}(x)$, $x \in \R{d}$,
\begin{equation}
  \label{eq:Harnack}
  \sup_{y \in B_1(x)}u(y) \leq c_H \inf_{y \in B_1(x)}u(y)\,,
\end{equation}
which shows that 
\begin{equation}
\label{eq:Harnack-1}
\mathbb{E}\left[\sup{_{x \in B_{1}(0)}P_{x,\omega}[\tilde{T}_{-\frac{L}{2}}^{l} < T_{\frac{L}{3}}^{l}]}\right] \leq c_H~ P_{0}[\tilde{T}_{-\frac{L}{2}}^{l} < 
T_{\frac{L}{3}}^{l}].
\end{equation}
Inserting (\ref{eq:Harnack-1}) in $(\ref{eq:P})$, we see that (\ref{eq:transience}) follows from (\ref{eq:geometric-Tgamma}).
From an application of Borel-Cantelli's lemma we obtain that $P_{0}$-a.s. 
for large integer L,
\[
T^{l}_{\frac{4L}{3}} < \tilde T_{\frac{L}{2}}^{l}
\circ \theta _{T_{L}^{l}} + T_{L}^{l}.
\]
So on a set of full $P_{0}$-measure we can construct an integer-valued 
sequence $L_{k} \nearrow \infty$, with $L_{k+1}=[\frac{4}{3}L_{k}]$ and
$T^{l}_{L_{k+1}} < \tilde{T}_{\frac{L_{k}}{2}}^{l}
\circ \theta _{T_{L_{k}}^{l}} + T_{L_{k}}^{l}, \, k \geq 0.$
This shows (\ref{eq:transience-1}).\\
We now show that for some $\mu>0$
\begin{equation}
\label{assertion}
\hat{E_{0}}[\exp\{{\mu \sup_{0 \leq t \leq \tau_{1}}|\Xt|^{\gamma}\}}] < \infty.
\end{equation}
The proof is divided into several propositions. In a first step, we study
the integrability properties of the random variable (recall (\ref{eq:return}))
\begin{equation}
\label{eq:M} 
M \df \sup{\{(\Xt-X_{0}) \cdot l : 0 \leq t \leq J\}}\, , 
\end{equation}
i.e. M is the maximal relative displacement of $X_{.}$ in the direction l
before it goes an amount of R below its starting point. 
By virtue of the Proposition \ref{prop:equiv} and (\ref{eq:transience-1}), we know that $P_{0}[D=\infty]=P_{0}[J=\infty]>0$ 
(recall (\ref{eq:return})).
Hence we cannot expect M to be finite. Nevertheless we have the following Proposition:
\begin{prop}
\label{prop:maximum}
There is $ \mu_1> 0$ such that 
\[
E_{0}[\exp{\{\mu_{1}~M^\gamma\}},~J < \infty] \leq 
1- \tfrac{P_{0}[J=\infty]}{2}.
\]
\end{prop}
\begin{proof} 
Let $L_k = \left(\frac{4}{3}\right)^k$. By our previous result (\ref{eq:transience}), we see that there is $\mu>0$ such that for large integers $k$:
\begin{equation}
\label{eq:M-0}
P_{0}\left[L_k \leq M < L_{k+1},~ J < \infty\right]  
\leq P_{0}\left[\tilde{T}^{l}_{L_k/2} \circ 
        \theta _{T_{L_k}^{l}} < 
        T^{l}_{4L_k/3} \circ \theta_{T_{L_k}^{l}}\right] 
\leq \exp{\{-\mu L_k^\gamma\}}.
\end{equation}
Let $k_{0}$ be large enough such that $\sum_{k \geq k_{0}} \exp{\{-\tfrac{\mu}
{2}L_k^\gamma\}} \leq \tfrac{P_{0}[J=\infty]}{4}.$ Further, let 
 $\mu_{1} >0$ such that $0 < (\frac{4}{3})^\gamma \mu_{1} < \frac{\mu}{2}$.
Then (\ref{eq:M-0}) shows that for $k_0$ large enough, 
\begin{eqnarray*}
&&E_{0}[\exp{\{\mu_1~M^\gamma\}},~J < \infty] \\ 
&\leq& \exp{\{\mu_1L_{k_0}^\gamma\}}~ P_{0}[J < \infty] +
  \sum_{k \geq k_{0}} \exp{\{\mu_1 L_{k+1}^\gamma\}}~
   P_{0}[L_k \leq M < L_{k+1},~ J < \infty] \\
&\leq&  \exp{\{\mu_1L_{k_0}^\gamma\}}~(1-P_{0}[J = \infty])+\sum_{k \geq k_{0}} \exp{\{-\tfrac{\mu}{2}L_k^\gamma\}}\\
&\leq& \exp{\{\mu_1L_{k_0}^\gamma\}}~(1-P_{0}[J = \infty])+ \tfrac{P_{0}[J=\infty]}{4} \leq 1-\tfrac{P_{0}[J=\infty]}{2},
\end{eqnarray*} 
provided  $\mu_1 >0$ is chosen small enough in the last inequality.
\end{proof}
As a next step, we shall prove the integrability of $\exp{\{\mu~ (X_{\tau_{1}}
\cdot l)^\gamma\}}$ under the extended annealed measure $\hat{P}_0$.
Recall the $(\mathcal{Z}_{t})_{t \geq 0}$- stopping times 
$(V_{k}(a))_{k \geq 0}$, $(\tilde{N}_{k}(a))_{k \geq 0}$ and
$N_{1}(a)$ defined in subsection \ref{subsec:regeneration}. As we will see in the proof of Proposition \ref{prop:Xtau1},
$\exp{\{\mu~((X_{N_{1}(a)}-X_{0}) \cdot l)^\gamma\}}$ will play a key role in studying the integrability of \\ $\exp{\{\mu~(X_{\tau_{1}} \cdot l)^\gamma\}}$ under
$\hat{P}_0$. Let us therefore start with the following Proposition, which
 only assumes that, $P_0-a.s.,\,\lim_{t \to \infty}X_t \cdot l =\infty$, which we have established in (\ref{eq:transience-1}).
\begin{prop}
\label{prop:N1}
Assume that $\lim_{t \to \infty}X_t \cdot l =\infty \,\, P_0-$a.s.
Then, for each $\mu_2 > 0$ there is $\mu_3 > 0$, such that
for $\mathbb{P}$-a.e. $\omega \in \Omega$:
\begin{equation}
\label{eq:N1}
\sup_{x,a \geq R}{\hat{E}_{x, \omega}}[\exp{\{\mu_3~(((X_{N_{1}(a)}-X_{0}) \cdot l)^\gamma-a^\gamma)\}}] \leq 1+\mu_2.
\end{equation}
\end {prop}
\begin{proof} 
Define $A_l \df \{\lim_{t \to \infty}\Xt \cdot l = \infty\}$. Observe that
\begin{equation}
\label{eq:t}
\text{for $\mathbb{P}$-a.e. $\omega$ and for every $x \in \R{d}$, 
$P_{x,\omega}[A_l]=1$}\,.
\end{equation}
Indeed, by the stationarity of the measure $\mathbb{P}$, $P_{y}[A_l]=1$ for all  $y \in \R{d}$. 
Hence $\int{dy~P_{y}[A_l^{c}]}=0$, and by applying Fubini's Theorem it follows that there is a $\mathbb{P}$-null set $\Gamma \subset \Omega$, such that for all
$\omega \notin \Gamma$ and y outside a Lebesgue null set $\mathcal{N}(\omega)
\subset \R{d}$, $P_{y,\omega}[A_l^{c}]=0$. 
Observe that for all $x \in \R{d}$, and  $\om \in \Omega$, $P_{x,\omega}[A_l]=P_{x,\omega}[A_l \circ \theta_{1}]$. 
It follows from the Markov property that for all $x \in \R{d}$, and $\om \notin \Gamma$,
$P_{x,\omega}[A_l \circ \theta_{1}]=\int_{\R{d}}P_{y,\omega}[A_l]\,p_\om(1,x,y)\,dy=1$,
where $p_\om(s,x,y)$ is the transition density function under $P_{x,\omega}$ (that is, for every open subset $U$ of $\R{d}$, $P_{x,\omega}[X_s \in U]=\int_U p_\om(s,x,y)\,dy$). The claim (\ref{eq:t}) now follows.
When $P_{x,\om}[A_l]=1$ for all $\om \in \Omega$ and all $x \in \R{d}$, Proposition 4.8 in Shen \cite{shen} shows that (\ref{eq:N1}) holds for all $\om \in \Omega$, when $\gamma =1$. By the same proof as given there, Proposition \ref{prop:N1} follows from (\ref{eq:t}) when $\gamma =1$.
When $0<\gamma<1$, using $\beta^\gamma-\alpha^\gamma
\leq \beta -\alpha$ for $\beta \geq 1 \vee \alpha$, and (\ref{eq:N1}) with 
$\gamma =1$, we find $\mu_3 \in (0,1)$ such that
\begin{align*}
  &\sup_{x,a \geq R}{\hat{E}_{x, \omega}}[\exp{\{\mu_3~(((X_{N_{1}(a)}-X_{0}) 
  \cdot l)^\gamma-a^\gamma)\}}]\\
  \leq & \sup_{x,a \geq R}{\hat{E}_{x, \omega}}[\exp{\{\mu_3~(((X_{N_{1}(a)}-X_  {0}) \cdot l)^\gamma-a^\gamma)\}},(X_{N_{1}(a)}-X_{0}) \cdot l \geq 1 \vee     a] + e \leq 4.
\end{align*}
By Jensen's inequality, if $n \geq 1$ is large enough, we find
\begin{align*}
     \sup_{x,a \geq R}{\hat{E}_{x, \omega}}[\exp{\{\tfrac{\mu_3}{n}~
     (((X_{N_{1}(a)}-X_{0}) \cdot l)^\gamma-a^\gamma)\}}]
\leq 4^{\frac{1}{n}} \leq 1+\mu_2\,,
\end{align*}
which shows (\ref{eq:N1}).
\end{proof}

\begin{prop}
\label{prop:Xtau1}
There exists $\mu_4 > 0$ such that 
\begin{equation}
\label{eq:tau-1}
\hat{E_{0}}[\exp\{{\mu_4 (X_{\tau_{1}} \cdot l)^{\gamma}\}}] < \infty.
\end{equation}
\end{prop}
\begin{proof} 
Using that, $\hat P_0$-a.s., $X_{S_{k}} \cdot l \leq X_{N_{k}}\cdot l+10R$, $k \geq 1$, (see the remark following (\ref{eq:N-S-R}))
we observe that 
\begin{multline}
\label{eq:h-k}
\hat E_{0}[\exp\{\mu_4 (X_{\tau_{1}} \cdot l)^\gamma\}] 
= \sum_{k \geq 1} \hat E_{0}[\exp\{{\mu_4 (X_{S_{k}} \cdot l)^\gamma\}},~S_{k}<\infty, D \circ \theta_{S_{k}}=\infty] \\ 
\leq \exp(\mu_4(10R)^\gamma) \sum_{k \geq 1} \hat E_{0}[\exp\{{\mu_4 (X_{N_{k}} \cdot l)^\gamma\}}, ~N_{k}<\infty]  
\df \exp(\mu_4(10R)^\gamma)\sum_{k \geq 1} h_{k}.
\end{multline}
Observe that, for $k \geq 1$, see (\ref{eq:N-S-R}),
\begin{equation*}
  l\cdot X_{N_{k+1}} = l\cdot X_{R_{k}}+l \cdot (X_{N_1(a_k)}-X_0)\circ\theta_{R_{k}}\,,
\end{equation*}
with $a_k=M(R_{k})-l\cdot X_{R_{k}}+R\in \mathcal Z_{R_k}$, (in fact for any $m\geq 1$, $a_k\cdot 1_{\{R_k=m\}}$ is $\mathcal F_m\otimes \mathcal S_{m-1}$-measurable, and $\lambda_m$ is independent of $\mathcal F_m\otimes\mathcal S_{m-1}$). We recall that the shift $\theta_{R_k}$ is {\em not} applied to $a_k$. Therefore, by the strong Markov property, cf. Proposition \ref{prop:coupling},
and, by applying Proposition \ref{prop:N1} (notice that $a_k \geq R$, $k \geq 1$, see (\ref{eq:N-S-R})), we see that for all $\mu_2>0$, there is $\mu_4 \in
(0,\mu_3)$ such that:
\begin{equation}
  \label{eq:integral-1}
  \begin{split}
  h_{k+1}
  \leq & \mathbb E\left[\hat E_{0,\om}\Big[\exp\big\{\mu_4(l\cdot X_{R_{k}})^\gamma\big\}, \, R_{k}<\infty,\,\hat E_{X_{R_{k}},\om}\big[\exp\big\{\mu_4(l\cdot (X_{N_1(a_k)}-X_0))^\gamma\big\}\big]\Big]\right]\\
  \leq & \mathbb E\left[\hat E_{0,\om}\Big[\exp\big\{\mu_4(l\cdot X_{R_{k}})^\gamma\big\}, \, R_{k}<\infty,\,(1+\mu_2)\,e^{\mu_4a_k^\gamma}\Big]\right]\,.
  \end{split}
\end{equation}
Observe that with $M$ from (\ref{eq:M}) and $Z_1$ as in Lemma \ref{lemma:bernstein} of the Appendix, the following inequalities hold, when $R_k$ is finite: 
\begin{gather*}
  a_k \leq Z_1\circ\theta_{J}\circ\theta_{S_k}+M\circ\theta_{S_{k}}+2R\,,\\
  l\cdot X_{R_{k}} = l\cdot X_{S_k}+\underbrace{\big(l\cdot (X_D-X_0)\big)}_
  {\leq Z_1\circ\theta_{J}}\circ~\theta_{S_k}\,.
\end{gather*}
Insert them into the last term of (\ref{eq:integral-1}), apply the strong Markov property at time $S_k$, cf. Proposition \ref{prop:coupling}, (we use the same argument as above, that for $m\geq 1$, $\exp\{\mu_4(l\cdot X_{S_k})^\gamma\}\cdot 1_{\{S_k=m\}}$ is $\mathcal F_m\otimes\mathcal S_{m-1}$-measurable, and $\lambda_m$ is independent of $\mathcal F_m\otimes\mathcal S_{m-1}$),  then use the strong Markov property for the process $(X_t)_{t\geq 0}$ at time $J$ on the event it is finite, and obtain (observe that $M$ is $\mathcal F_{J}$-measurable)
\begin{align*}
  &  h_{k+1}    \\
  \leq& e^{\mu_4(2R)^\gamma}\, \mathbb E\left[\hat E_{0,\om}\Big[e^{\mu_4(l\cdot X_{S_k})^\gamma}, S_k<\infty, (1+ \mu_2)\, \hat E_{X_{S_k},\om}\big[\exp\big\{\mu_4(2Z_1^\gamma \circ\theta_J+M^\gamma)\big\}, J<\infty\big]\Big]\right]\\
  \leq& e^{\mu_4(2R)^\gamma}\,\mathbb E\left[\hat E_{0,\om}\Big[e^{\mu_4(l\cdot X_{S_k})^\gamma}, S_k<\infty, (1+ \mu_2)\, E_{X_{S_k},\om}\Big[e^{\mu_4M^\gamma}\, E_{X_J,\om}\big[e^{2\mu_4Z_1^\gamma}\big], J<\infty\Big]\Big]\right]\,.
\end{align*}  
From Lemma \ref{lemma:bernstein} of the Appendix, we know that, for $\mu_4 \in(0,\delta)$, $\sup_{x,\omega} E_{x,\om}[e^{2\mu_4Z_1^\gamma}]\leq 1+ \mu_2$. 
Further we use that,  $\hat P_0$-a.s.,$(X_{S_{k}}-X_{N_{k}})\cdot l \leq 10R$, $k \geq 1$, and, that, after an application of the strong Markov property to the stopping time $N_k$, conditionally on $\mathcal Z_{N_k}$, $X_1$ is uniformly 
distributed on $B^{X_{N_k}}$ under $\hat P_{X_{N_k},\om}^1$, see Proposition 
\ref{prop:coupling}. Let $\mu_5 \df \exp\{\mu_4((2R)^\gamma+(10R)^\gamma)\}
(1+ \mu_2)^2$, 
 then we obtain that the last expression is smaller than
\begin{align*}
&\mu_5\,  \mathbb E\left[\hat E_{0,\om}\Big[e^{\mu_4(l\cdot X_{N_k})^\gamma}, 
N_k<\infty, E_{X_{S_k},\om}\Big[e^{\mu_4M^\gamma}\,,J<\infty\Big]\Big]\right]\\
=& \mu_5 \, \mathbb E\left[\hat E_{0,\om}\Big[e^{\mu_4(l\cdot X_{N_k})^\gamma}, 
N_k<\infty, \hat E_{X_{N_k},\om}^1 \Big[E_{X_{1},\om}\Big[e^{\mu_4M^\gamma}\,,
J<\infty \Big]\Big]\Big]\right]\\
=&\mu_5\frac{1}{|B_R|}\,\int \textrm{d}y \,
\mathbb E\left[\hat E_{0,\om}\Big[e^{\mu_4(l\cdot X_{N_k})^\gamma}, 
N_k<\infty, y \in B^{X_{N_k}}\Big] E_{y,\om}\Big[e^{\mu_4M^\gamma}\,,J<\infty
\Big]\right]\,.
\end{align*}
Since $\hat{E}_{0,\omega}[\exp{\{\mu_4 (X_{N_{k}} \cdot l)^\gamma\}},
~N_{k}<\infty,~y \in B^{X_{N_{k}}}]$ is $\mathcal{H}_{\{x \cdot l \leq 
y \cdot l-4R\}}$-measurable (see point (3) in the addendum \cite{shen-add}  to 
Shen \cite{shen}) and 
${E}_{y,\omega}[\exp{\{\mu_4~M^\gamma\}},~J < \infty]$ is $\mathcal{H}_{\{x 
\cdot l \geq y \cdot l-R\}}$-measurable, as a result of finite range dependence, see (\ref{eq:R-separation}), the above random variables are $\mathbb{P}$-independent.
Hence, using the stationarity of the measure $\mathbb{P}$ and Proposition
\ref{prop:maximum}, we obtain that
\begin{align*}
h_{k+1} 
\leq & \mu_5 \hat E_0\Big[e^{\mu_4(l\cdot X_{N_k})^\gamma}, N_k<\infty\Big]
\cdot E_0\Big[e^{\mu_4M^\gamma}\,,J<\infty\Big]\\
\leq& \mu_5 \Big(1-\tfrac{P_{0}[J=\infty]}{2}\Big)\hat E_0\Big[e^{\mu_4(l\cdot X_{N_k})^\gamma}, N_k<\infty\Big]
\leq (1-\alpha) \hat E_0\Big[e^{\mu_4(l\cdot X_{N_k})^\gamma}, N_k<\infty\Big],
\end{align*}
for some $\alpha > 0$, provided $\mu_2 >0$ and $\mu_4 \in (0,\mu_1 \wedge \mu_3
\wedge \delta)$ are small enough such that 
$\mu_5(1-\tfrac{P_{0}[J=\infty]}{2})=e^{\mu_4((2R)^\gamma+(10R)^\gamma)}(1+ \mu_2)^2(1-\tfrac{P_{0}[J=\infty]}{2}) \leq 1-\alpha$.
It follows by induction that
\[
h_{k+1} \leq (1-\alpha)^{k}\hat E_0\Big[e^{\mu_4(l\cdot X_{N_1})^\gamma}, N_1<\infty\Big]
\]
so that, by (\ref{eq:h-k}), and by virtue of Proposition \ref{prop:N1},
\[
\hat{E_{0}}\Big[e^{\mu_4 (X_{\tau_{1}} \cdot l)^\gamma}\Big] \leq e^{\mu_4(10R)^\gamma}
\hat E_0\Big[e^{\mu_4(l\cdot X_{N_1})^\gamma}, N_1<\infty\Big]~\sum_{k \geq 0} (1-\alpha)^{k} < \infty\,,
\]
which is our claim (\ref{eq:tau-1}).
\end{proof}
The assertion (\ref{assertion}) now readily follows. Choose $r>0$ such that
$ \mathcal{ \bar D} \subset B_{r}(0)$, and let $\tilde{L} \df \frac{L^{\frac{1}{\gamma}}}{r}$. Hence $\tilde{L}\mathcal{ \bar D} \subset B_{r\tilde{L}}(0)$, and by definition of the random variable $\tau_1$ in (\ref{eq:regeneration-time}),
\begin{equation}
\label{eq:sup}
\begin{aligned}
 &\hat{P}_{0}\left[\sup_{0 \leq t \leq \tau_{1}}{|\Xt|^{\gamma}} \geq L\right]
\leq \hat{P}_{0}[T_{\tilde{L}\mathcal{D}} < \tau_{1}]
\leq P_0[T_{\tilde{L}\mathcal{D}} < 
                   T^{l}_{\tilde{L}}]
      + \hat{P}_{0}[T_{\tilde{L}\mathcal{D}} = 
                   T^{l}_{\tilde{L}},~
                 T_{\tilde{L}\mathcal{D}} < \tau_{1} ]\\
\leq& P_0[T_{\tilde{L}\mathcal{D}} < T^{l}_{\tilde{L}}] +
      \hat{P}_{0}[X_{\tau_{1}} \cdot l \geq \tilde{L}-3R].
\end{aligned}
\end{equation}
Applying (\ref{eq:exit-D}) to the first term on the right-hand side of
(\ref{eq:sup}), and applying Chebychev's inequality and Proposition \ref{prop:Xtau1} to the second term on the right-hand side, we find
\begin{align}
\label{eq:assertion-1}
  \limsup_{L \to \infty} L^{-1} \log \hat P_0[\sup_{0 \leq t \leq 
\tau_{1}}{|\Xt|^{\gamma}} \geq L]<0.
\end{align}
Thus, for some $\mu>0$ small enough,
\begin{equation*}
\hat{E}_{0}[\exp{\{\mu \sup_{0 \leq t \leq \tau_{1}}|\Xt|^{\gamma} \}}]
= 1+\mu \int_{0}^{\infty}\exp{\{\mu L\}}\hat{P}_{0}[\sup_{0 \leq t \leq 
\tau_{1}}{|\Xt|^{\gamma}} \geq L]~dL < \infty\,,
\end{equation*}
and (\ref{assertion}) follows from (\ref{eq:assertion-1}).

\subsection{The Proof of (\ref{eq:integrability-2}) $\Rightarrow$ (\ref{eq:geometric-Tgamma})}
By Proposition \ref{prop:equiv}, we know that $\lim_{t \to \infty} \Xt \cdot l
= \infty$ $P_{0}$-a.s. implies $\tau_{1} < \infty$ $\hat{P}_{0}$-a.s., and hence Theorem \ref{thm:renewal} holds.
To verify condition $(T)_{\gamma}|l$, we first show that the diffusion has an asymptotic direction $\hat{v}$ under $\hat{P}_{0}$, with $\hat{v} \cdot l >0$, see
Proposition \ref{prop:asymptotic-direction}. The claim (\ref{eq:geometric-Tgamma}) is implied by Lemma \ref{lemma:exit}, which is immediate for $d=1$, and which follows from a control on the oscillations of the diffusion orthogonal to $\hat{v}$ under $\hat{P}_{0}$, see Proposition \ref{prop:orthogonal}, when $d \geq 2$.
\begin{prop}
\label{prop:asymptotic-direction}
It holds  that
\begin{equation}
\label{eq:asymptotic-direction}
P_0-\text{a.s.},\, \,\,\frac{\Xt}{|\Xt|} \underset{t \to \infty}{\longrightarrow} \hat{v} \df
\frac{\hat{E}_{0}[X_{\tau_{1}}|D=\infty]}{|\hat{E}_{0}[X_{\tau_{1}}|D=\infty]|}
\qquad \textrm{and} \qquad \hat{v} \cdot l >0.
\end{equation}
\end{prop}
\begin{proof} 
By definition of $\tau_{1}$, $X_{\tau_{1}} \cdot l >0$ $\hat{P}_{0}$-a.s., so $\hat{v}$ is well defined and $\hat{v} \cdot l >0.$ 
By assumption, $\hat{E}_{0}[X_{\tau_{1}}|D=\infty]<\infty$. The strong
law of large numbers applied to the i.i.d. random variables
$X_{\tau_{k+1}}-X_{\tau_{k}}$, $k \geq 1$, (cf. Theorem \ref{thm:renewal})
 yields
\begin{equation}
\label{eq:law-LN}
\frac{1}{k}X_{\tau_{k}} \underset{k \to \infty}{\longrightarrow}
\hat{E}_{0}[X_{\tau_{1}}|D=\infty] \qquad \hat{P}_{0}-\textrm{a.s.}
\end{equation}
For $t>0$, define $k(t)$ via 
\begin{equation}
\label{eq:number-of-regenerations}
\tau_{k(t)} \leq t < \tau_{k(t)+1},
\end{equation}
i.e. $k(t)$ is the number of regenerations up to time t.
Clearly $\hat P_0$-a.s. $k(t)\underset{t \to \infty}{\longrightarrow} \infty$. Write, for $k(t) \geq 1$,
\begin{equation}
\label{eq:remainder-term}
\frac{\Xt}{k(t)}=\frac{X_{\tau_{k(t)}}}{k(t)}+\frac{1}{k(t)}
(\Xt-X_{\tau_{k(t)}}).
\end{equation}
The modulus of the second term on the right-hand side can be bounded by
\begin{equation}
 \sup_{s \geq 0}\frac{1}{k(t)}\, \big| X_{(\tau_{k(t)}+s) \wedge \tau_{k(t)+1}}-X_{\tau_{k(t)}}\big|.
\end{equation}
Since $\lambda_{\tau_{k}-1}=1$, $k \ge 1$, it follows from Proposition 2.1 that, $\hat P_0$-a.s., $X_u \in U^{X_{\tau_{k}-1}}$ for all 
$u \in [\tau_k-1,\tau_k]$, and we thus find that $\hat P_0$-a.s.,
\begin{equation}
\label{path decomp}
\tfrac{1}{k}\, \big| X_{(\tau_{k}+s) \wedge \tau_{k+1}}-X_{\tau_{k}}\big|
\leq  \tfrac{1}{k}\, \big|X_{(\tau_{k}+s) \wedge (\tau_{k+1}-1)}-X_{\tau_{k}}\big| + \tfrac{12R}{k}\,.
\end{equation}
For $k \geq 0$, let $Y_k \df \sup_{s \geq 0} |X_{(\tau_{k}+s) \wedge (\tau_{k+1}-1)}-X_{\tau_{k}}|$.
From Theorem \ref{thm:renewal}, we know that the random variables $Y_{k},\, 
k \geq 1$, are i.i.d. random variables under $\hat{P}_{0}$ and are distributed under $\hat P_0$ as $Y_{0}$ under $\hat P_0[\cdot|D=\infty]$.
Hence, applying Chebychev's inequality and Theorem \ref{thm:renewal}, we find
by virtue of (\ref{eq:integrability-2}) that, for $\epsilon >0$, there is
$\mu >0$ and $\alpha < \infty$ such that  for $k \geq 1$,
\begin{multline*}
\hat{P}_{0}\Big[\tfrac{|Y_{k}|}{k}>\epsilon \Big]
\leq \exp{\{-\mu (k\epsilon)^{\gamma}\}}\hat{E}_{0}
[\exp{\{\mu|Y_{k}|^{\gamma}\}}]\\
=\exp{\{-\mu (k\epsilon)^{\gamma}\}}\hat{E}_{0}
[\exp{\{\mu \sup_{s \geq 0}|X_{s \wedge (\tau_{1}-1)}|^{\gamma}\}}\, \big| D=\infty]
\leq \alpha \, \exp{\{-\mu (k\epsilon)^{\gamma}\}}.
\end{multline*}
Applying Borel-Cantelli's lemma, we see that, $\hat{P}_{0}$-a.s., 
$\frac{1}{k}|Y_{k}|\underset{k \to \infty}{\longrightarrow} 0$, and hence,
$\hat{P}_{0}$-a.s.,
$\frac{1}{k(t)}|Y_{k(t)}| \underset{t \to \infty}{\longrightarrow} 0$. The claim (\ref{eq:asymptotic-direction}) now follows from 
(\ref{eq:law-LN}), (\ref{eq:remainder-term})and from (\ref{path decomp}).
\end{proof}
Denote by $\Pi(\, \cdot \,)$ the orthogonal projection on the orthogonal complement of $\hat{v}$:
\begin{equation}
\label{eq:orth-proj}
 \Pi(w) \df w-(w \cdot \hat{v})\hat{v}\,,
\end{equation}
and let $L_{u}^{l} \df \sup \{t \geq 0: \Xt \cdot l \leq u\}$ be the time 
of last visit of the half space $\{x \cdot l \leq u\}$ by $X_\cdot$. 
The next Proposition gives a control on the oscillations of the process
orthogonal to $\hat{v}$, when $d \geq 2$.
\begin{prop}
\label{prop:orthogonal}
($d \geq 2$)\\
Assume (\ref{eq:integrability-2}). For $\rho \in (\frac{1}{2},1]$ and $\alpha>0$,
\begin{equation}
\label{orthogonal}
\limsup_{u \to \infty} u^{-(2\rho-1)\wedge \gamma \rho} \log P_{0}
\Big[\sup_{0 \leq t \leq L_{u}^l}|\Pi(\Xt)|>\alpha\, u^{\rho}\Big]<0.
\end{equation}
\end{prop}
\begin{proof}
Without loss of generality, we can replace $|\Pi(\Xt)|$ by $\Xt \cdot w$,
where $w \in \R{d}$ is such that $w \cdot \hat{v}=0$. Recall the definition
of $k(t)$ in (\ref{eq:number-of-regenerations}). 
Notice that $\hat{P}_{0}$-a.s.,
 for $k \geq 1$, $(X_{\tau_{k}}-X_{\tau_{k-1}})\cdot l \geq 21R/2$.
Indeed, by Theorem \ref{thm:renewal}, it suffices to prove the statement for $k=1$. Recall that $\tau_0=0$, and observe that $(X_{V_k (3R)}-X_0)\cdot l \ge 3R$, all $k \ge 0$, and hence we find $(X_{\tilde N_1 (3R)}-X_0)\cdot l \ge 5R/2$, all $k\ge 1$.
Consequently, since $X_{\tilde N_k (3R)} \cdot l \ge X_{\tilde N_1 (3R)}\cdot l$, we obtain
$(X_{N_1 (3R)}-X_0)\cdot l \ge 5R/2$, and since $\lambda_{N_1 (3R)}=1$, we find
from Proposition \ref{prop:coupling}, as well as 
from the definition of $\tau_1$ and the stopping times $S_k$, $k \ge 1$, that 
$(X_{\tau_1}-X_0)\cdot l \ge (X_{S_1}-X_0)\cdot l \ge (X_{N_1(3R)}-X_0)\cdot l +8R \ge 21R/2$.
Since, for $0 \leq t \leq L^l_u$, $X_{\tau_{k(t)}} \cdot l < u+R$, it follows that
$k(t) \leq \frac{u+R}{21R/2} \leq \frac{u}{R}$, u large enough.
Let $X^\ast \df \sup_{t \leq \tau_1}|X_t-X_0|$. For $t \geq 0$ 
it holds $\hat P_0$-a.s. that
\begin{equation*}
  X_t \cdot w=X_{\tau_{k(t)}}\cdot w+(X_t-X_{\tau_{k(t)}})\cdot w
  \leq X_{\tau_{k(t)}}\cdot w+X^\ast \circ \theta_{\tau_{k(t)}}.
\end{equation*}
It follows that 
\begin{equation}
\label{eq:X*}
\begin{aligned}
 & \hat P_0[\sup_{0 \leq t \leq L_{u}}X_t \cdot w > \alpha u^\rho]
\leq \sum_{0 \leq k \leq \frac{u}{R}} \hat P_0[X_{\tau_k}\cdot w+X^\ast \circ \theta_{\tau_k}> \alpha u^\rho]\\
\leq & \sum_{0 \leq k \leq \frac{u}{R}} \hat P_0[X^\ast \circ \theta_{\tau_k}> \tfrac{\alpha}{3} u^\rho]+\sum_{1 \leq k \leq \frac{u}{R}}(\hat P_0[X_{\tau_1}
\cdot w >\tfrac{\alpha}{3} u^\rho]+\hat P_0[(X_{\tau_k} - X_{\tau_1})\cdot w
>\tfrac{\alpha}{3} u^\rho]).
\end{aligned}
\end{equation}
Applying first Chebychev's inequality, then Theorem \ref{thm:renewal} to the first term of the last line of (\ref{eq:X*})(we use the same decomposition of the path as in (\ref{path decomp})), and with (\ref{eq:integrability-2}) applied 
 to both the first and the second term, we find that there is $\lambda > 0$,
 such that for large $u$, (\ref{eq:X*}) is smaller than 
\begin{align}
\label{eq:sum-1}
  \exp\{-\lambda (\tfrac{\alpha}{3}u^{\rho})^\gamma\} + \sum_{1 \leq k \leq \frac{u}{R}}
  \hat P_0[(X_{\tau_k} - X_{\tau_1})\cdot w >\tfrac{\alpha}{3} u^\rho].
\end{align}
If $\gamma \in (0,1)$, the claim (\ref{orthogonal}) follows from Theorem 
\ref{thm:renewal} and from Theorem A.1.
in the Appendix of Sznitman \cite{szn02}. If $\gamma=1$, then, as above,
we first apply Chebychev's inequality and then Theorem \ref{thm:renewal} to (\ref{eq:sum-1})
and obtain that it is smaller than
\begin{equation*}
  \exp\{-\lambda \tfrac{\alpha}{3}u^\rho\} (1+ \sum_{1 \leq k \leq \frac{u}{R}}
  \hat E_0 [\exp\{\lambda X_{\tau_1}\cdot w\}|D=\infty]^{k-1})
\leq \exp\{-\lambda \tfrac{\alpha}{3}u^\rho\}(1+\tfrac{u}{R} \exp\{\tfrac{u}{R}H(
\lambda)\}),
\end{equation*}
provided, we define, for $|\lambda|$ small, 
\begin{equation*}
  H(\lambda) \df \log \,\hat E_0 [\exp\{\lambda X_{\tau_1}\cdot w\}|D=\infty].
\end{equation*}
$H(\cdot)$ is a convex function, and, since $\hat E_0[X_{\tau_1}\cdot w\,|\,D=\infty]=0$, we see that $H(0)=0,\, 
H'(0)=0,\,H(\cdot) \geq 0$ for $\lambda \geq 0$, and $H(\lambda)=O(\lambda^2)$, as $\lambda \to 0$. If $\rho =1$, choose $\lambda >0$ small enough such that
$H(\lambda)<\lambda \tfrac{\alpha}{3} R$, and (\ref{orthogonal}) holds. In the case $\rho \in (\tfrac{1}{2},1)$, we instead choose for a sufficiently small 
$\nu >0,\, \lambda=\nu u^{\rho -1}$, and conclude in a similar fashion.
\end{proof}
Let $\hat{R}(\cdot)$ be a rotation of $\R{d}$ such that $\hat{R}(e_{1})=\hat v$. For $\epsilon > 0$, consider the cylinder in $\R{d}$:
\begin{equation}
\label{rotation}
C^{\epsilon, u} \df \hat{R}\left(\left(-\epsilon u,\frac{u}{\epsilon}\right)
\times B^{d-1}_{\frac{\epsilon u}{2}}(0)\right),
\end{equation}
where, for $r>0$, $B^{d-1}_r(0)$ stands for the $(d-1)$-dimensional Euclidean ball with radius $r$ and center 0.
($C^{\epsilon, u}$ is understood as $\hat{R}(-\epsilon u,\frac{u}{\epsilon})$
when $d=1$).\\
The next step is
\begin{lemma}
\label{lemma:exit}
Assume (\ref{eq:integrability-2}). For $\epsilon > 0$, 
\begin{equation}
\label{eq:exit-1}
\limsup_{u \rightarrow \infty} u^{-\gamma} \log P_{0}[T_{C^{\epsilon, u}} <
T^{\hat{v}}_{\frac{u}{\epsilon}}]<0.
\end{equation}
\end{lemma}
\begin{proof}
Let us first handle the case $d=1$. From Chebychev's inequality and (\ref{eq:geometric-Tgamma}), for large $u$, we find $\alpha >0$ such that
\begin{equation*}
 P_0[\widetilde T^{\hat{v}}_{-u \epsilon}<T^{\hat{v}}_{\frac{u}{\epsilon}}]
\leq  P_0[\widetilde T^{\hat{v}}_{-u \epsilon} < \infty] \leq
\hat P_0 [\sup_{0 \leq t \leq \tau_1}|X_t| \geq \epsilon u] \leq \exp\{-\alpha u^\gamma\}\,,
\end{equation*}
and (\ref{eq:exit-1}) follows.
 When $d \geq 2$, write 
\begin{multline}
\label{eq:sum-2}
  P_0[T_{C^{\epsilon, u}} <T^{\hat{v}}_{\frac{u}{\epsilon}}] \leq
  P_0[\widetilde T^{\hat{v}}_{-u \epsilon}<T^{\hat{v}}_{\frac{u}{\epsilon}},\,
  \sup \{|\Pi(X_t)|:t \leq \widetilde T^{\hat{v}}_{-u \epsilon}\} \leq
  \tfrac{\epsilon}{2}\,l \cdot \hat v \,u] +\\
  P_0[\widetilde T^{\hat{v}}_{-u \epsilon}<T^{\hat{v}}_{\frac{u}{\epsilon}},\,
   \tfrac{\epsilon}{2}\,l \cdot \hat v \,u < \sup \{|\Pi(X_t)|:t \leq 
    \widetilde T^{\hat{v}}_{-u \epsilon}\} \leq \tfrac{\epsilon}{2}u]+
  P_0[T_{C^{\epsilon, u}} <\widetilde T^{\hat{v}}_{-u \epsilon} \wedge
  T^{\hat{v}}_{\frac{u}{\epsilon}}].
\end{multline}
Let us first estimate the probability of the leftmost event on the right-hand side of (\ref{eq:sum-2}). Observe that on this event,
\begin{equation*}
  X_{\widetilde T^{\hat{v}}_{- \epsilon u}}\cdot l=
  X_{\widetilde T^{\hat{v}}_{- \epsilon u}}\cdot \hat v\,\hat v \cdot l
+\Pi(X_{\widetilde T^{\hat{v}}_{- \epsilon u}})\cdot l
\leq -\frac{\epsilon}{2}u \,\hat v \cdot l\,.
\end{equation*}
Hence, with the help of (\ref{eq:integrability-2}), for large $u$, we find $\alpha >0$ such that the probability of this event is smaller than 
\begin{equation*}
P_0 [\widetilde T^{l}_{-\tfrac{\epsilon}{2}l \cdot \hat v \,u}<\infty]
\leq   \hat P_0 [\tau_1 > \widetilde T^{l}_{-\tfrac{\epsilon}{2}l \cdot 
 \hat v \,u}] \leq \hat P_0 [\sup_{0 \leq t \leq \tau_1}|X_t| \geq 
  \tfrac{\epsilon}{2}l \cdot \hat v \,u]
 \leq \exp\{-\alpha u^\gamma\}.
\end{equation*}
To bound the rightmost term of (\ref{eq:sum-2}), notice that
$\{T_{C^{\epsilon, u}}<\widetilde T^{\hat{v}}_{-u \epsilon} \wedge
  T^{\hat{v}}_{\frac{u}{\epsilon}}\} \subseteq $ \\
$\{\sup_{0 \leq t \leq L^l
_{(\epsilon/2 +1/\epsilon)u}}|\Pi(X_t)|\geq \tfrac{\epsilon u}{2}\}$, 
and then apply Proposition \ref{prop:orthogonal} with $\rho =1$. The bound
for the middle term of (\ref{eq:sum-2}) equally follows from a direct application of Proposition \ref{prop:orthogonal} with $\rho =1$.
\end{proof}
Now (\ref{eq:geometric-Tgamma}) easily follows. Indeed, choose
$\epsilon > 0$ such that $\epsilon < 2b \wedge \frac{\hat{v}\cdot l}{2}$. 
The last estimate also holds for unit vectors $l'$ in a neighborhood of $l$, and, with 
 the notation $\partial_+  C^{\epsilon, L}=\{x \in \partial \,C^{\epsilon, L}: x \cdot \hat v=L/\epsilon\}$ for the ``top part'' of the boundary of the cylinder
and similarly $\partial_-  C^{\epsilon, L}=\{x \in \partial \,C^{\epsilon, L}: x \cdot \hat v=-\epsilon L\}$ for the ``bottom part'' of the boundary, it follows
that 
$\partial_+  C^{\epsilon, L}$ is contained in the complement of $U_{l',b,L}$, whereas 
$\partial_-  C^{\epsilon, L}$ lies inside $U_{l',b,L}$.
As a result, we find that for unit vectors $l'$ as above,
\begin{equation}
\label{cylinder}
\limsup_{L \rightarrow \infty} L^{-\gamma} \log P_{0}\left[X_{T_{U_{l',b,L}}}
\cdot l' < 0\right] 
\leq \limsup_{L \rightarrow \infty} L^{-\gamma} \log P_{0}\left[
T_{C^{\epsilon, L}} < T^{\hat{v}}_{\frac{L}{\epsilon}}\right]
<0\,,
\end{equation}
which is our claim (\ref{eq:geometric-Tgamma}).

\begin{rem}
\label{direction T}  
In the same way as in (\ref{cylinder}), we see that,  
\begin{equation}
\label{eq:T-half-plane}
\text{if $(T)_{\gamma}|l_0$ holds
for some $l_0 \in S^{d-1}$, then $(T)_{\gamma}|l$ holds 
 iff $l \cdot \hat{v}>0$}\,.
\end{equation}
\end{rem}

\section{Tail estimates on the first renewal time $\tau_{1}$}
\label{sec:tail-estimate}
The ballistic law of large numbers and the central limit theorem established
in Shen \cite{shen} (see (\ref{eq:lln}) and (\ref{eq:clt})) respectively follow from
 $P_0$-a.s. $\lim_{t \to \infty} X_t \cdot l=\infty$, 
$\hat E_0[\tau_1] < \infty$ and from $P_0$-a.s. $\lim_{t \to \infty} X_t \cdot l=\infty$, $\hat E_0[\tau_1^2] < \infty$.
In this section, we are going to derive tail estimates on $\tau_{1}$ under the assumption of condition $(T')$. These will  ensure the finiteness of every moment
of $\tau_{1}$ when $d \geq 2$, see (\ref{eq:tail-estimate}).
The arguments in this section closely follow section 3 in Sznitman \cite{szn01}.\\
For a bounded domain $U$, and $f$ a bounded measurable function on $U$, introduce the semigroup corresponding to the diffusion killed when exiting $U$, see (\ref{eq:exit-time}) for notations, 
\begin{equation}
\label{eq:semigroup}
R^{U}_{t,\om}f(x) \df E_{x,\om}[f(X_{t}), T_{U}>t],
\end{equation}
and a threshold time related to the decay of the semigroup,
\begin{equation}
\label{treshhold time}
t_{\om}(U) \df \inf\left\{t \geq 0:\|R^{U}_{t,\om}\|_{\infty,\infty}\leq \frac{1}{2}\right\}
=\inf \left\{t \geq 0:\sup_{x \in U}P_{x,\om}[T_U>t]\leq \frac{1}{2}\right\}.
\end{equation}
Consider further the successive returns of $X_{\cdot}$ to $B_{1}(x)$ and 
departures from $B_{2}(x)$,  
\begin{equation}
\label{eq:excursion-1}
 R_{1}^{x} \df \inf \{s \geq 0: X_{s} \in B_{1}(x)\},~~ 
 D_{1}^{x} \df \inf \{s \geq  R_{1}^{x}: X_{s} \notin B_{2}(x)\},
\end{equation}
and inductively, for $n \geq 0$,
\begin{equation}
\label{eq:excursion-2}
 R_{n+1}^{x} \df D_{n}^{x}+R_{1}^{x}\circ \theta_{D_{n}^{x}},~~
 D_{n+1}^{x} \df R_{n+1}^{x}+D_{1}^{x}\circ \theta_{R_{n+1}^{x}}.
\end{equation}
\begin{lemma}
\label{lemma:threshold-time}
There is a constant $c$ such that for all  bounded domains $U$ and $\om \in \Omega$, one can find $x_{0}$ in $\frac{1}{\sqrt{d}}\mathbb{Z}^{d}$ within distance 1 of $U$ such that
\begin{equation}
\label{eq:threshold-time}
\inf_{z \in \partial B_2 (x_0)}P_{z,\om}[R_{1}^{x_{0}}>T_{U}] \leq \frac{c \, \text{diam}(U)^d}{t_{\om}(U)}\,.
\end{equation}
\end{lemma}
\begin{proof}
Cover $U$ by unit balls centered in $\tfrac{1}{\sqrt{d}}\mathbb Z^d$, and let 
$(y_{i})_{i=1}^{N}$, $N \leq c\,\text{diam}(U)^d$, be an enumeration of the centers of these balls. 
Choose $\delta \leq t_\om (U)/2$, then, by definition of $t_{\om}(U)$, we can find
an $x_{1}$ in $U$ such that $P_{x_{1},\om}[T_{U}>t_{\om}(U)-\delta]>\frac{1}{2}$. Hence $\frac{1}{4}t_{\om}(U) \leq \frac{1}{2}(t_{\om}(U)-\delta) \leq E_{x_{1},\om}[T_{U}]$.
Applying the strong Markov property to the stopping times $R_{j}^{y_i}$ 
and using the fact that $\sup_{\om \in \Omega} \sup_{i,x \in \bar B_{1}(y_{i})}
E_{x,\om}\left[T_{B_{2}(y_{i})}\right]<\infty$, see for instance \cite{kar-shr} p.365, yields
\begin{multline}
\label{eq:return-1}
\frac{1}{4}t_{\om}(U) \leq  E_{x_{1},\om}[T_{U}]
\leq \sum_{i=1}^{N} E_{x_{1},\om}\left[\int_{0}^{T_{U}}\mathbf{1}_{B_{1}(y_{i})}
(X_{s})\textrm{d}s\right]\\
\leq  \sum_{i=1}^{N}\sum_{j=1}^{\infty}E_{x_{1},\om}\left[R_{j}^{y_{i}}<T_{U},
E_{X_{R_{j}^{y_{i}}},\om}\left[\int_{0}^{D_{1}}\mathbf{1}_{B_{1}(y_{i})}
(X_{s})\textrm{d}s\right]\right]
\leq \,c\,\sum_{i=1}^{N}\sum_{j=1}^{\infty}P_{x_{1},\om}\left[R_{j}^{y_{i}}<T_{U}\right].
\end{multline}
For $j \geq 2$, successive applications of the strong Markov property show that
\begin{align*}
  P_{x_{1},\om}[R_{j}^{y_{i}}<T_{U}]=&
E_{x_1, \om}[R_{j-1}^{y_{i}}<T_{U}, \,P_{X_{D^{y_i}_{j-1}},\om}[R_1^{y_i}<T_U]]\\
\leq & \sup_{z \in \partial B_2(y_i)} P_{z,\om} [R_1^{y_i}<T_U]\,P_{x_1,\om}
[R_{j-1}^{y_{i}}<T_{U}] \\
\leq &(\sup_{z \in \partial B_2(y_i)} P_{z,\om} 
[R_1^{y_i}<T_U])^{j-1}P_{x_1,\om}[R_{1}^{y_{i}}<T_{U}]\,.
\end{align*}
Using the last estimate, we see that the last expression in (\ref{eq:return-1}) is smaller than
\begin{align*}
c~\sum_{i=1}^{N} \frac{P_{x_{1},\om} [R_{1}^{y_{i}}<T_{U}]}
{\inf_{z \in \partial{B_{2}(y_{i})}}P_{z,\om} [R_{1}^{y_{i}}>T_{U}]}
\leq \frac{c\,\text{diam($U$)}^{d}}{\inf_{1 \leq i \leq N}
\inf_{z \in \partial{B_{2}(y_{i})}}P_{z,\om} [R_{1}^{y_{i}}>T_{U}]}\,.
\end{align*}
The claim (\ref{eq:threshold-time}) now follows.
\end{proof}

For $\beta \in (0,1]$ and $L>0$, we denote by $U_{\beta,L}$ the set
\[
U_{\beta,L} \df \{x \in \R{d}: x \cdot l \in (-L^{\beta},L)\}.
\]
The next Proposition shows that the control of the tail of the variable
$\tau_{1}$ can be obtained from the derivation of large-deviation-type estimates on the exit distribution of the diffusion out of $U_{\beta,L}$.
\begin{prop}
\label{prop:trap}
Let $d \geq 2$, and assume that $(T')$ holds with respect to $l \in S^{d-1}$.
If $\beta \in (0,1)$ is such that for any $\alpha>0$,
\begin{equation}
\label{eq:trap}
\limsup_{L \rightarrow \infty} L^{-1}\log \mathbb{P}\left[P_{0,\om}\left[
X_{T_{U_{\beta,L}}}\cdot l >0\right]\leq \exp\{-\alpha L^{\beta}\}\right]<0,
\end{equation}
then
\begin{equation}
\label{eq:tail-2}
\limsup_{u \rightarrow \infty}\,(\log u)^{-\zeta}\log \hat P_{0}[\tau_{1}>u]<0 
\end{equation}
for any $\zeta<\frac{1}{\beta}$ (when $(T)$ holds, one can choose
$\zeta=\frac{1}{\beta}$).
\end{prop}
\begin{proof}
Let $R$ be a rotation of $\R{d}$ such that $R(e_{1})=l$. For $L>0$ write
\begin{eqnarray*}
C_{L}=R\left(\left(-L/2,L/2\right)^{d}\right) 
\quad \mathrm{and}\quad
V_{x}=x+R\left((-1,3)\times (-1,1)^{d-1}\right).
\end{eqnarray*}
From the Support Theorem, see \cite{bass} p.25, 
we know that there is a constant
$\kappa >0$ such that for all $x \in \R{d}$ and all $\om \in \Omega$
\begin{equation}
\label{kappa}
\inf_{z \in B_{\frac{1}{2}}(x)}P_{z,\om}[X_{1} \in B_{\frac{1}{2}}(x+2l),
T_{V_{x}}>1]\geq \kappa >0.
\end{equation}
For $u>1$, denote 
$\Delta(u) \df \lfloor \frac{\log u}{6 \log(1/\kappa)}\rfloor \quad \mathrm{and} \quad
L(u) \df \Delta(u)^{\frac{1}{\beta}}$.
Let $\beta \in (0,1)$ and $\zeta < \frac{1}{\beta}$. Write
\begin{equation}
\label{eq:tail-tau-1}
\begin{aligned}
\hat{P}_{0}[\tau_{1}>u] \leq & \hat{P}_{0}[\tau_{1}>u,~T_{C_{L(u)}} \leq \tau_{1}]
+P_{0}[T_{C_{L(u)}}>u]\\
\leq &\hat{P}_{0}[\sup_{0 \leq t \leq \tau_{1}}|X_t| \geq L(u)/2]+
P_{0}[T_{C_{L(u)}}>u].
\end{aligned}
\end{equation}
Using Chebychev's inequality and condition $(T)_\gamma|l$, $\gamma$ close to 1 
such that $\frac{\gamma}{\beta}\geq \zeta$, we find that
\begin{equation}
  \label{eq:first-term}
  \limsup_{u \rightarrow \infty}\, (\log u)^{-\zeta} \log 
   \hat{P}_{0}[\sup_{0 \leq t \leq \tau_{1}}|X_t| \geq L(u)/2]<0.
\end{equation}
Hence, by means of (\ref{eq:tail-tau-1}), it suffices to show that 
\begin{equation}
\label{eq:tail-1}
\limsup_{u \rightarrow \infty}\left(\log u\right)^{-\frac{1}{\beta}}
\log P_{0}[T_{C_{L(u)}}>u]<0.
\end{equation}
Recall the definition of $t_{\om}(U)$ in (\ref{treshhold time}), and denote by
$\mathcal{T}$ the event
\begin{equation}
\mathcal{T} \df \left\{\om \in \Omega:t_{\om}(C_{L(u)})>\frac{u}{(\log u)^{\frac{1}{\beta}}}\right\}.
\end{equation}
It follows from Lemma \ref{lemma:threshold-time} and the 
Markov property that for large $u$
\begin{align}
\label{eq:tau}
&P_{0}[T_{C_{L(u)}}>u] \leq \mathbb{E}\left[\mathcal{T}^{c},
P_{0,\om}[T_{C_{L(u)}}>u]\right]+\mathbb{P}[\mathcal{T}]\nonumber \leq \\
& \left(\frac{1}{2}\right)^{\lfloor (\log u)^{\frac{1}{\beta}}\rfloor}+
\mathbb{P}\Big[\exists~ x_{2} \in C_{L(u)} \cap \tfrac{1}{\sqrt{d}}\mathbb{Z}^{d};
\inf_{z \in \partial B_2(x_2)}P_{z,\om}[R_{1}^{x_{2}}>T_{C_{L(u)}}]
 \leq \frac{c L(u)^d (\log u)^{\frac{1}{\beta}}}{u}\Big].
\end{align}
(Notice that, if $x_2$ would not belong to $C_{L(u)}$, then we would find
from the Support Theorem, see \cite{bass} p.25, that for every $z \in \partial
B_2(x_2)$, $P_{z,\om}[R_{1}^{x_{2}}>T_{C_{L(u)}}]\geq c >0$, which contradicts
the rightmost event in the last line for large $u$.) 
Choose $x=x_{2}+2 \Delta(u)l$.
By the strong Markov property, we see that
\begin{equation}
\label{eq:escape-prob}
\inf_{z \in \partial B_2(x_2)}P_{z,\om}[R_{1}^{x_{2}}>T_{C_{L(u)}}]
\geq \inf_{z \in \partial B_2(x_2)}P_{z,\om}[R_{1}^{x_{2}}>R_{1}^{x}]~
\inf_{z \in \partial B_{1}(x)}P_{z,\om}[R_{1}^{x_{2}}>T_{C_{L(u)}}].
\end{equation}
Let $y \in \partial B_{2}(x_{2})$. One way to hit $B_{1}(x)$ before returning to $B_{1}(x_{2})$ when starting at $y$ is the following: we hit
$B_{\frac{1}{2}}(x_2+2l)$ before hitting $B_{1}(x_{2})$ which happens with probability at least $\tilde \kappa$, where $\tilde \kappa$ is a positive constant, see the Support Theorem p.25 in \cite{bass}. Then we hit 
$B_{\frac{1}{2}}(x_{2}+4l)$ without exiting
$V_{x_2+2l}$ which occurs with probability at least $\kappa$, see (\ref{kappa}). Then continue hitting $B_{\frac{1}{2}}(x_{2}+2(k+1)l)$ without exiting
$V_{x_{2}+2kl}$, $1 \leq k \leq \Delta(u)-1$,
until landing in $B_{1}(x)$. Hence 
\begin{equation}
\label{kappa2}
\inf_{z \in \partial B_2(x_2)}P_{z,\om}[R_{1}^{x_{2}}>R_{1}^{x}] \geq 
\tilde{\kappa} \kappa^{\Delta(u)-1}
\geq \tilde \kappa u^{-\frac{1}{6}}.
\end{equation}
Together with (\ref{eq:escape-prob}), this shows that for large $u$, on the event $\mathcal T$, see (\ref{eq:tau}),
\begin{eqnarray}
\label{eq:exit-before-return}
\inf_{z \in \partial B_{1}(x)}P_{z,\om}[R_{1}^{x_{2}}>T_{C_{L(u)}}]\leq 
\frac{1}{\tilde \kappa}u^{\frac{1}{6}}
\inf_{z \in \partial B_2(x_2)}P_{z,\om}[R_{1}^{x_{2}}>T_{C_{L(u)}}]
\leq u^{-\frac{1}{2}}. 
\end{eqnarray}
In particular, by a similar argument as given below (\ref{eq:tau}), we see that, for large $u$,  $B_{3}(x) \subset C_{L(u)}$.
By the same argument as in (\ref{eq:strong-Markov}), it follows that, for large $u$,
$P_{\cdot,\om}[R_{1}^{x_{2}}>T_{C_{L(u)}}]$ is  
$\mathcal L_\om$-harmonic on $B_{3}(x)$, and (\ref{eq:Harnack}) shows that
\begin{equation}
\label{eq:harnack-1}
P_{x,\om}[R_{1}^{x_{2}}>T_{C_{L(u)}}] \leq c_H \inf_{z \in \partial B_{1}(x)}
P_{z,\om}[R_{1}^{x_{2}}>T_{C_{L(u)}}]\,.
\end{equation}
It follows from (\ref{eq:exit-before-return}) and (\ref{eq:harnack-1}) that
for large $u$,
\begin{eqnarray*}
P_{x,\om}[X_{T_{x+U_{\beta,L(u)}}}\cdot l>x \cdot l] \leq
P_{x,\om}[R_{1}^{x_{2}}>T_{C_{L(u)}}]
\leq c_H u^{-\frac{1}{2}}
\leq \exp\big(-cL(u)^{\beta}\big).
\end{eqnarray*}
Using translation invariance and (\ref{eq:tau}), we find
\[
P_{0}[T_{C_{L(u)}}>u] \leq \left(\frac{1}{2}\right)^{\lfloor (\log u)^{\frac{1}{\beta}}\rfloor}+ c\, L(u)^{d}\,\mathbb{P}\Big[P_{0,\om}[X_{T_{U_{\beta,L(u)}}}\cdot l>0]\leq \exp \big(-c L(u)^{\beta}\big)\Big],
\]
and (\ref{eq:tail-1}) follows from (\ref{eq:trap}). This proves (\ref{eq:tail-2}).
\end{proof}
We shall now derive upper bounds like (\ref{eq:trap}) under the assumption of condition $(T')$. By means of Proposition \ref{prop:trap}, we then obtain tail estimates on the first renewal time $\tau_1$.  
We first need some notation. For $\beta > 0$ and $L>0$, consider the lattice 
\begin{equation*}
\mathcal{L}_{\beta,L}=L\mathbb{Z}\times((2d+1)L^{\beta}+2R)\mathbb{Z}^{d-1},
\end{equation*}
and, for $w \in \R{d}$, we introduce the blocks
\begin{equation}
\label{blocks}
\begin{aligned}
&B_{1,\beta,L}(w)=\hat{R}(w+[0,L]\times[0,L^{\beta}]^{d-1}),\\
&B_{2,\beta,L}(w)=\hat{R}(w+(-dL^{\beta},L]\times(-dL^{\beta},(d+1)L^{\beta})
^{d-1}),
\end{aligned}
\end{equation}
where $\hat{R}$ is a rotation of $\R{d}$ such that $\hat{R}(e_1)=\hat v$, and $\hat v$
is the asymptotic direction of the annealed diffusion (that exists under $(T')$, see Proposition \ref{prop:asymptotic-direction}). We shall 
also consider the following subset of the boundary of $B_{2,\beta,L}(w)$,
which is a subset of the 'top part' of the box,
\[
\partial_{+}{B}_{2,\beta,L}(w)=\partial{B}_{2,\beta,L}(w) \cap \partial B_{1,\beta,L}(w), \quad w \in \R{d},
\]
as well as the random variables
\begin{equation}
X_{\beta, L}(w)=-\log \inf_{x \in B_{1,\beta,L}(w)}P_{x,\om}\left[
X_{T_{B_{2,\beta,L}(w)}} \in \partial_{+}{B}_{2,\beta,L}(w)\right].
\end{equation}
To obtain an upper bound like (\ref{eq:trap}) under $(T')$, it is instrumental to produce a control on the 
tail of the random variable $X_{\beta,L}(w)$ for some $\beta \in (0,1)$ under $(T')$.
Indeed, we devise an escape route for the diffusion through the ``right'' side of $U_{\beta,L}$ by
piling up in the direction $\hat v$ a finite number of boxes of type $ B_{2,\beta,L}$. 
An atypical behavior of the exit distribution out of the slab $U_{\beta,L}$ under $P_{0,\om}$ as in  (\ref{eq:trap})
implies an atypical size for at least one of the $X_{\beta,L}(w)$ in one of the piled up boxes. Hence, 
to produce an upper bound like (\ref{eq:trap}), it suffices to 
show that, for large $L$, the probability that $X_{\beta,L}(w)$ is bigger than const $L^\beta$ decays exponentially with $L$ for some $\beta \in (0,1)$.\\ 
We prove in fact a stronger statement. Namely, we show that the above probability decays exponentially with $L^\zeta$, where $\zeta < f(\beta)=d(2 \beta -1)$, with $\beta$ 
restricted to the interval $(1/2,1)$, so that for suitable values of $\beta$ close to 
one, $\zeta$ can be chosen larger than one, since $d \ge 2$.
By means of a renormalization-type argument, see Lemma \ref{lemma:renormalisation}, we reduce this task to showing a substantially weaker estimate. Indeed, it now suffices to prove for some $\beta_0$ slightly larger than $1/2$ that the probability that $X_{\beta_0,L}(w)$ is bigger than const $L^\beta$ decays exponentially with $L^{f_0(\beta)}$, where $f_0(\beta)=\beta+\beta_0 -1$, and $\beta \in (\beta_0,1)$. 
This ``seed-estimate'' is then provided in Lemma \ref{lemma:seed-estimate} under the assumption of condition $(T')$.\\
We begin with the renormalisation step. 
Surprisingly enough, we do not need to assume condition $(T')$, in which
case the rotation $\hat{R}$ in (\ref{blocks}) is an arbitrary rotation of $\R{d}$.
\begin{lemma}[Renormalisation step, $d \geq 2$]
\label{lemma:renormalisation}
Assume that $\beta_{0} \in (0,1)$ and $f_{0}$ is a positive function defined on $[\beta_{0},1)$, such that 
\begin{equation*}
f_{0}(\beta)\geq f_{0}(\beta_{0})+\beta - \beta_{0}, \quad \beta \in [\beta_{0},1)
\end{equation*}
and, for $\beta \in [\beta_{0},1)$, $\zeta < f_{0}(\beta)$,
\begin{equation}
\lim_{\beta' \uparrow \beta}\limsup_{L \rightarrow \infty}L^{-\zeta}
\sup_{w \in \R{d}} \log \mathbb{P}[X_{\beta_{0},L}(w)\geq L^{\beta'}]<0.
\end{equation}
Denote by $f(\cdot)$ the linear interpolation on $[\beta_{0},1]$ of the 
value $f_{0}(\beta_{0})$ at $\beta_{0}$ and the value $d$ at 1.Then, for
$\beta \in [\beta_{0},1)$ and $\zeta < f(\beta)$,
\begin{equation}
\label{eq:renormalisation}
\lim_{\beta' \uparrow \beta}\limsup_{L \rightarrow \infty}L^{-\zeta}
\sup_{w \in \R{d}} \log \mathbb{P}[X_{\beta,L}(w)\geq L^{\beta'}]<0.
\end{equation}
\end{lemma}
\begin{proof}
We only give a sketch of the proof, since it is similar to the proof of Lemma 3.2 in \cite{szn01}. For $\chi \in (0,1)$ defined via $\beta \df \chi \beta_0 + 1-\chi$, we consider the set
\begin{equation}
  \label{eq:Col}
  \text{Col} \df \{z \in \mathcal L_{\beta_0, L^\chi}, z \cdot e_1=0,
                   z \cdot e_i \in [\frac{1}{4}L^\beta,\frac{3}{4}L^\beta],\,
2 \leq i \leq d\}\,. 
\end{equation}
For $w \in \R{d}$, attach at every $w+z$, $z \in$ Col, a ``column of boxes''
$B_{1,\beta_0,L^\chi}(\cdot)$, made by piling up $\lfloor L^{1-\chi}\rfloor$
such boxes on top of each other. Each such column will provide a line of escape of the diffusion out of a box $B_{2,\beta,L}(w)$ through $\partial_+ B_{2,\beta,L}(w)$. Every $x \in B_{1,\beta,L}(w)$ is at most at distance $\sqrt{d}L^\beta$ from a box $B_{1,\beta_0,L^\chi}(\cdot)$ in one of the aforementioned
columns. From a similar argument as in (\ref{kappa2}), and from the strong Markov property, we see that for large $L$ and $c_1=\sqrt{d} \log\frac{1}{\kappa}$, with $\kappa$ from (\ref{kappa}), $J=\lfloor  L^{1-\chi} \rfloor$,
\begin{equation*}
  \{ X_{\beta,L}(w) \geq 3 c_1 L^\beta\} \subseteq \{\min_{z \in \text{Col}}
  \underbrace{\sum_{j=0}^{J} X_{\beta_0,L^\chi}(w+z  +jL^\chi e_1)}_{\df Y(z)} \geq 2 c_1 L^\beta\}\,.
\end{equation*}
Using the independence of the variables $Y(z)$, $z \in$ Col, and Chebychev's inequality, we find that for $\lambda >0$, 
\begin{equation}
\label{eq:indep-1}
  \mathbb P[X_{\beta,L}(w) \geq 3 c_1 L^\beta] \leq 
  \underset{z \in \text{Col}}{\Pi} \left\{\exp\{-\lambda c_1 L^\beta\} \mathbb E
   [\exp\{\tfrac{\lambda}{2}Y(z)\}]\right\}\,.
\end{equation}
Observe that, for $z \in$ Col and large $L$, the variables $X_{\beta_0,L^\chi}(w+z  +jL^\chi e_1)$ are independent when $j$ is restricted to the set of even or the set of odd integers. It thus follows from Cauchy-Schwarz's inequality 
that the right-hand side of (\ref{eq:indep-1}) is smaller than 
\begin{equation*}
\underset{z \in \text{Col}}{\Pi}\left\{\exp\{-\lambda c_1 L^\beta\}\Pi_{j=0}^{J} \mathbb E[\exp\{\lambda X_{\beta_0,L^\chi}(w+z  +jL^\chi e_1)\}]^{1/2} \right\}\,.
\end{equation*}
Since the random variables $X_{\beta_0,L^\chi}$ are non-negative, 
the quantity in the last line becomes larger when we omit the square roots,
and an application of Fubini's Theorem yields that the last line can be bounded by
\begin{equation}
\label{eq:indep-2}
\underset{z \in \text{Col}}{\Pi}\Big\{\exp\{-\lambda c_1 L^\beta\}
       \Big(\exp\{\tfrac{\lambda}{2} c_1 L^{\chi \beta_0}\}+
       \int_{\tfrac{c_1}{2}L^{\chi \beta_0}}^{\infty}\lambda
       e^{\lambda u} \sup_{w' \in \R{d}}\mathbb P[X_{\beta_0,L^\chi}(w') \geq u]\textrm{d}u\Big)^{J+1}\Big\}\,.
\end{equation}
For $\lambda=L^{\alpha}$, $\alpha = \chi f_0(\beta_0)-\chi \beta_0-\varepsilon$ and $0<\varepsilon < \chi f_0(\beta_0)$, one can show that the integral in the rightmost term of (\ref{eq:indep-2}) tends to 0 as $L \to \infty$. Since 
$\lambda L^{\chi \beta_0}$ tends to $\infty$ with $L$, we find that, for large $L$, 
\begin{equation*}
  \sup_{w \in \R{d}} \mathbb P[X_{\beta,L}(w) \geq 3 c_1 L^\beta] \leq
  \exp\{-\tfrac{\lambda}{6}c_1 L^\beta \text{\#Col}\}\,.
\end{equation*}
Since \#Col$\sim c L^{(d-1)(\beta-\chi \beta_0)}$, as $L \to \infty$, we
obtain that, for small $\varepsilon >0$,
\begin{equation}
  \label{eq:escape}
  \limsup_{L \to \infty} L^{-(\chi f_0(\beta_0)+d(1-\chi)-\varepsilon)}
  \sup_{w \in \R{d}}\log \mathbb P[X_{\beta,L}(w)\geq 3 c_1 L^\beta]<0\,,
\end{equation}
which implies the claim.
\end{proof}
The next Lemma shows that, when $d \geq 2$, under condition $(T')$, the
function $f_{0}(\beta)=\beta +\beta_{0}-1$, $\beta \in [\beta_{0},1)$,
fulfills the assumption of Lemma \ref{lemma:renormalisation} when $\beta_{0} \in (\frac{1}{2},1)$.
\begin{lemma}[Seed estimate, $d \geq 2$, under $(T')$]
\label{lemma:seed-estimate}
Assume that $\beta_{0} \in (\frac{1}{2},1)$. Then, for $\rho > 0$ and $\beta \in [\beta_{0},1)$,
\begin{equation}
\label{eq:seed-estimate}
\limsup_{L \rightarrow \infty}L^{-(\beta +\beta_{0}-1)}
\sup_{w \in \R{d}} \log \mathbb{P}[X_{\beta_{0},L}(w)\geq \rho L^{\beta}]<0.
\end{equation}
\end{lemma}
\begin{proof}
Choose $\eta \in (0,1)$ small and then introduce $\chi=\beta_{0}+1-\beta \in
(\beta_{0},1]$, and, for large $L$ and $w \in \R{d}$ the boxes $\tilde{B}_{1}(w) \subset \tilde{B}_{2}(w)$, defined analogously as before, with $[0,L]\times[0,L^{\beta}]^{d-1}$ and $(-dL^{\beta},L]\times(-dL^{\beta},(d+1)L^{\beta})
^{d-1})$ replaced by $[0,L_{0}]\times[0,L^{\beta_{0}}]^{d-1}$ and 
$(-dL^{\beta_0},L_{0}+3]\times(-\eta L^{\beta_{0}},(1+\eta )L^{\beta_{0}})^{d-1})$ 
respectively, with the notation
\[
L_{0}=\frac{L-\eta L^{\beta_{0}}}{\lfloor L^{1-\chi}\rfloor}.
\]
Define also Top$\tilde{B}_{2}(w)=\partial{\tilde{B}}_{2}(w) \cap \{x:
x \cdot \hat{v}=w \cdot \hat{v}+L_{0}+3\}$.
Let $(B_1(z_i))_{i \in I}, z_i \in \tilde B_1 (w)$, 
$I$ a finite set
growing polynomially with $L$, be a finite cover of $\tilde B_1 (w)$ by unit balls. For $L$ large,
it holds that $B_3(z_i) \subset \tilde B_2(w), i \in I$, and 
by the same argument as in (\ref{eq:strong-Markov}), we see that
$P_{\cdot,\om}\big[ X_{T_{\tilde{B}_{2}(w)}}
\in \mathrm{Top}~\tilde{B}_{2}(w)\big]$ is $\mathcal L_\om$-harmonic on
$B_3(z_i)$, so that (\ref{eq:Harnack}) implies that for all $i \in I$,
\begin{align}
\label{eq:harnack-2}
P_{z_i,\om}\left[ X_{T_{\tilde{B}_{2}(w)}}\in \mathrm{Top}~\tilde{B}_{2}(w)
\right] \leq c_H \inf_{x \in  B_{1}(z_i)}P_{x,\om}\left[ X_{T_{\tilde{B}_{2}(w)}}
\in \mathrm{Top}~\tilde{B}_{2}(w)\right].
\end{align}
We say that $w$ is {\it good} when
\[
\inf_{x \in \tilde{B}_{1}(w)}P_{x,\om}\left[ X_{T_{\tilde{B}_{2}(w)}}
\in \mathrm{Top}~\tilde{B}_{2}(w)\right]\geq \frac{1}{2\,c_H},
\]
and {\it bad} otherwise.
Hence, by (\ref{eq:harnack-2}), and using Chebychev's inequality and translation invariance, we obtain
\begin{equation}
\label{eq:bad}
\begin{aligned}
\mathbb{P}[w \text{ is bad}] 
\leq & \sum_{i \in I} \mathbb{P}\Big[\inf_{x \in B_{1}(z_i)}P_{x,\om}
\big[ X_{T_{\tilde{B}_{2}(w)}}\in \mathrm{Top}~\tilde{B}_{2}(w)\big]
<\frac{1}{2\,c_H}\Big] \\
\leq &\sum_{i \in I}\mathbb{P}\left[P_{z_i,\om}\left[ X_{T_{\tilde{B}_{2}(w)}}
\in \mathrm{Top}~\tilde{B}_{2}(w)\right]<\frac{1}{2}\right] \\
\leq & 4|I|\Big(P_{0}[\sup_{0 \leq t \leq T^{\hat{v}}_{L_{0}+3}}|\Pi(X_{t})|\geq \eta L^{\beta_{0}}]+
P_{0}[\widetilde{T}^{\hat{v}}_{-dL^{\beta_{0}}}<\infty]\Big).
\end{aligned} 
\end{equation}
Notice that $L_0 \sim L^\chi$, so that, for large $L$, $T^{\hat{v}}_{L_{0}+3} \le T^{\hat{v}}_{2L^\chi+3}
\le  L^{\hat{v}}_{2L^\chi+3}$.
Then, under condition $(T)_{\gamma}|l$,
where $\gamma$ fulfills $\gamma \beta_{0} \geq 2\beta_{0}-\chi$,
we find with the help of Proposition \ref{prop:orthogonal} applied (with $\rho=\beta_0/\chi \in (1/2,1)$ and 
$u=2L^\chi +3$)
to the first term on the right-hand side of (\ref{eq:bad}) that
\begin{equation}
\label{eq:orthogonal}
\limsup_{L \rightarrow \infty} L^{-(2\beta_{0}-\chi)}\log P_{0}
[\sup_{0 \leq t \leq T^{\hat{v}}_{L_{0}+3}}|\Pi(X_{t})|\geq \eta L^{\beta_{0}}]
<0,
\end{equation}
and, since $(T)_{\gamma}|\hat v$ holds, see (\ref{eq:T-half-plane}), we find 
with the help of Chebychev's inequality that there is $\mu >0$ such that
\begin{equation}
\label{eq:left}
P_{0}[\tilde{T}^{\hat{v}}_{-dL^{\beta_{0}}}<\infty] \leq \hat{P}_{0}
\Big[\sup_{0 \leq t \leq \tau_{1}}|X_{t}|\geq d L^{\beta_{0}}\Big]
\leq \exp (-\mu L^{\gamma \beta_{0}})\,.
\end{equation}
According to our choice of $\gamma$, we obtain with 
(\ref{eq:bad}),(\ref{eq:orthogonal}) and (\ref{eq:left}) that
\begin{equation}
\label{eq:bad-0}
\limsup_{L \rightarrow \infty} L^{-(2\beta_{0}-\chi)}\sup_{w \in \R{d}}
\log \mathbb{P}[w \mathrm{~~is~~bad}]<0.
\end{equation}
When starting in $B_{1,\beta_0,L}(w) \cap \tilde B_1(w+j_0L_0 e_1)$, $0 \leq j_0
 <\lfloor L^{1-\chi}\rfloor$, for large $L$, one way to exit $B_{2,\beta_0,L}(w)$ through $\partial_+ B_{2,\beta_0,L}(w)$ is to successively exit the boxes 
$\tilde B_2(w+jL_0e_1)$, $j_0 \leq j <\lfloor L^{1-\chi}\rfloor$, through 
Top $\tilde B_2(w+jL_0e_1)$, and move to the box  $\tilde B_1(w+(j+1)L_0 e_1)$, which is at distance at most $\sqrt d \,\eta L^{\beta_0}$ from every point in 
Top $\tilde B_2(w+jL_0e_1)$, until landing in $\tilde B_1(w+\lfloor L^{1-\chi}\rfloor L_0 e_1) \cap B_{1,\beta_0,L}(w)$, and then exit $B_{2,\beta_0,L}(w)$ through $\partial_+ B_{2,\beta_0,L}(w)$, which is at distance at most $\eta L^{\beta_0}$ from every point in $\tilde B_1(w+\lfloor L^{1-\chi}\rfloor L_0 e_1) \cap B_{1,\beta_0,L}(w)$. When $w \in \R{d}$ and all $w+jL_0e_1$,  $0 \leq j <\lfloor L^{1-\chi}\rfloor$, are good, then, for large $L$, it follows from the strong Markov property and from (\ref{kappa}) that for all $x \in B_{1,\beta_0,L}(w)$,
\begin{equation}
  P_{x,\om}[X_{T_{B_{2,\beta_0,L}(w)}} \in \partial_+ B_{2,\beta_0,L}(w)]
  \geq \Big( \frac{1}{2c_H}\kappa ^{\lceil \tfrac{1}{2}\sqrt d \, \eta L^{\beta_0}\rceil
  +1}\Big)^{L^{1-\chi}} \kappa^{\lceil \tfrac{\eta}{2}L^{\beta_0}\rceil
  +1}
  >\exp\{-\rho L^\beta\}\,,
\end{equation}
provided $\eta >0$ is chosen small enough such that $\tfrac{\eta}{2}(1+\sqrt d) \log\tfrac{1}{\kappa}<\tfrac{\rho}{2}$, where $\rho>0$ is as in (\ref{eq:seed-estimate}). Therefore, for large $L$,
\begin{equation*}
  \sup_{w \in \R{d}} \mathbb P[X_{\beta_0, L}(w) \geq \rho L^\beta]
\leq L^{1-\chi}  \sup_{w \in \R{d}}\mathbb P[w \text{ is bad}]\,,
\end{equation*}
and the claim (\ref{eq:seed-estimate}) follows from (\ref{eq:bad-0}) together
with the identity $2 \beta_0 - \chi= \beta_0 + \beta -1$.
\end{proof}
We can now state the main result. With the help of the Renormalisation Lemma \ref{lemma:renormalisation},
we propagate the seed estimate contained in Lemma \ref{lemma:seed-estimate} to the right scale, 
and by piling up a finite number of boxes of the type $B_{2,\beta,L}$ in the direction $\hat v$, we obtain an upper bound like (\ref{eq:trap}).
Proposition \ref{prop:trap} then enables us to obtain tail estimates on $\tau_1$.
\begin{thm}($d \geq 2$)
\label{thm:tail-estimate}
Assume that $(T')$ holds relative to $l$. Then, for $\beta \in (\frac{1}{2},1)$,
\begin{equation}
\label{eq:trap-1}
\limsup_{L \rightarrow \infty} L^{-\zeta} \log \mathbb{P}\,
[P_{0,\om}[X_{T_{U_{\beta, L}}}\cdot l >0]\leq \exp \{-L^{\beta}\}]<0
\mathrm{~~for}~~ \zeta < d(2\beta -1),
\end{equation}
and
\begin{equation}
\label{eq:tail-estimate}
\limsup_{u \rightarrow \infty}\,(\log u)^{-\alpha} \log \hat P_{0}[\tau_{1}>u]<0 
\mathrm{~~for~~} \alpha<1+\frac{d-1}{d+1}\,.
\end{equation}
\end{thm}
\begin{proof}
Let $\beta$ and $\zeta$ be as in (\ref{eq:trap-1}), and choose $\beta_0 \in (\frac{1}{2},\beta)$ close to $\tfrac{1}{2}$, as well as $\beta' \in (\beta_0,\beta)$ such that, in the notation of Lemma \ref{lemma:renormalisation}, $f(\beta') > \zeta$. By piling up $N$ boxes
$B_{1,\beta',L}, B_{2,\beta',L}$, $0 \leq j \leq N$, where $N$ is chosen as the smallest integer such that 
\begin{equation*}
  Nl\cdot \hat v >1\,,
\end{equation*}
we obtain from the strong Markov property that for large $L$, 
\begin{center}
$P_{0,\om}[X_{T_{U_{\beta, L}}}\cdot l >0] \geq \exp \left\{-\sum_{j=0}
^{N}X_{\beta',L}(jLe_1)\right\}$  so that\\
$\mathbb{P}\,[P_{0,\om}[X_{T_{U_{\beta, L}}}\cdot l >0]\leq \exp \{-L^{\beta}\}]\leq (N+1)\, \underset{w}{\sup}\, \mathbb{P}\left[X_{\beta',L}(w) \geq \frac{L^\beta}{N}\right]$\,.
\end{center}
(\ref{eq:trap-1}) now follows from (\ref{eq:renormalisation}) applied with
$f_0(\cdot)=\beta_0+\cdot-1$, in view of Lemma \ref{lemma:seed-estimate}. 
For the proof of (\ref{eq:tail-estimate}), let $\alpha \in (1,2d/(d+1))$, and define $\beta = \alpha^{-1}$. Then, for any $\mu>0$, 
\begin{equation*}
 \limsup_{L \rightarrow \infty} L^{-1} \log \mathbb{P}\,
 [P_{0,\om}[X_{T_{U_{\beta, L}}}\cdot l >0]\leq \exp \{-\mu L^{\beta}\}]<0\,,
\end{equation*}
as follows from (\ref{eq:trap-1}) applied to $\beta' \in (\tfrac{1}{2},\beta)$, such that $d(2 \beta' -1) >1$. The claim now follows from Proposition \ref{prop:trap}.
\end{proof}

\section{Examples of condition (T)}
\label{sec:examples}
We start with an easy example. 
\begin{prop}
\label{prop:non-nestling}
($d \geq 1$) If for some $\delta >0$ and all $\om \in \Omega$, 
all $x \in \R{d}$,
\begin{equation}
  \label{eq:non-nestling-0}
  b(x,\om)\cdot l > \delta\,,
\end{equation}
then condition $(T)|l$ holds.
\end{prop}
\begin{proof}
  Define for $u \in \R{}$, $s(u) \df \exp\{-\tfrac{\delta}{\nu}u\}$. It follows from (\ref{eq:b-sigma-bound}), (\ref{eq:non-nestling-0}) that $s(X_t \cdot l)$ is a supermartingale, and an
application of Chebychev's inequality and of the stopping theorem yield that for all $\om \in \Omega$
\begin{equation}
  P_{0,\om}[X_{T_{U_{l,b,L}}}\cdot l <0] \leq \frac{1}{s(-bL)}E_{0,\om}[
s(X_{T_{U_{l,b,L}}}\cdot l)]
\leq \exp\{-\tfrac{\delta b}{\nu}L\}\,.
\end{equation}
The set of unit vectors that satisy (\ref{eq:non-nestling-0}) is open, and hence condition $(T)|l$ holds.
\end{proof}
Consequently, when $d \geq 2$, we recover and extend the main result of Komorowski and Krupa \cite{kom-krupa}, which provides a law of large numbers when $\sigma = Id$.
Proposition \ref{prop:non-nestling} holds for a general diffusion matrix $\sigma$ that satisfies (\ref{eq:b-sigma-bound})-(\ref{eq:R-separation}),
and we have in addition a central limit theorem, see (\ref{eq:lln}) and (\ref{eq:clt}).\\

We will now turn to a more involved situation. In the remainder of this section we now assume that, cf. (\ref{eq:b-sigma-bound}), (\ref{eq:SDE}), (\ref{eq:diff-operator}),
\begin{equation}
  \label{eq:sigma}
  \sigma (\cdot) = Id \,.
\end{equation}
The next Theorem provides a rich class of examples of diffusions in random
environment which fulfill condition $(T)$, and hence, when $d \geq 2$, a ballistic law of large numbers, and a central limit theorem with non-degenerate 
covariance matrix governing corrections to the law of large numbers, see 
(\ref{eq:lln}) and (\ref{eq:clt}).
\begin{thm}
\label{thm:examples}
($d \geq 1$) Assume (\ref{eq:stationarity})-(\ref{eq:R-separation})
and (\ref{eq:sigma}). There is a constant $c_e>0$, such that for $l \in S^{d-1}$, 
\begin{equation}
\label{eq:criterion}
\E[(b(0,\omega) \cdot l)_{+}]>c_e~ \E[(b(0,\omega) \cdot l)_{-}]
\end{equation}
implies $(T)|l$, cf. (\ref{eq:T}).
\end{thm}
Theorem \ref{thm:examples} is the main result of this section. Its analogue
in the discrete i.i.d. setting can be found in \cite{bolt-szn} p.40. 
In contrast to Proposition \ref{prop:non-nestling}, it comprises situations
where $b(0,\om)\cdot l$ changes sign for every unit vector $l$, see also
remark \ref{rem:examples} at the end of this section.\\

The proof of Theorem \ref{thm:examples} is inspired by the strategy used
in the discrete i.i.d. setting, see \cite{bolt-szn} p.40. Following
Kalikow's idea, for each bounded domain $U$, we introduce an auxiliary diffusion with characteristics independent of the environment, see (\ref{eq:b'}) and (\ref{eq:measure}). When starting at 0, this diffusion and the annealed diffusion have the same exit
distribution from $U$, see Proposition \ref{prop:exit}. This restores some Markovian character to the question of controlling exit distributions of $X_\cdot$ under the annealed measure, and enables us to show that condition $(T)$ is
implied by a certain condition $(K)$, see (\ref{eq:K}), which has a similar
flavor as Kalikow's condition in the discrete i.i.d. setting, see \cite{szn-zer}. The proof of Theorem \ref{thm:examples} is then carried out by checking condition ($K$).\\
Let us now define the auxiliary diffusion process mentioned above. Let $U$ be a bounded domain containing 0, and, for $x,y \in U$, $s>0$, denote with $p_{\omega, U}(s,x,y)$ the subtransition density for the quenched diffusion started in $x$ and killed when exiting $U$ ($p_{\omega, U}(s,x,y)$ can for instance be defined by means of Duhamel's formula, see equation
(\ref{eq:inf-green}) in the appendix or \cite{stroock} page 331).
We define the corresponding Green function through
\begin{equation}
  \label{eq:green-function}
  g_U(x,y,\om) \df \int_{0}^{\infty} p_{\omega, U}(s,x,y)ds\,.
\end{equation}
We now define the auxiliary drift term 
\begin{equation}
\label{eq:b'}
b'_U(x) \df 
 \begin{cases}
  \frac{\E[g_{U}(0,x,\omega)b(x,\omega)]}{\E[g_{U}(0,x,\omega)]}\,,
  &\text{if $x \in U \smallsetminus \{0\}$}\,, \\
  0\,,  & \text{if $x=0$ or $x \in U^c$\,.}
 \end{cases}
\end{equation}
The next lemma will be useful in the sequel. 
\begin{lemma}
\label{lemma:b'}
It holds that $|b'_U(x)|\leq \bar b$ (see (\ref{eq:b-sigma-bound}) for the notation), and $g_U(0,\cdot,\om)$ and $b'_U(\cdot)$ are continuous in $U \smallsetminus \{0\}$.
\end{lemma}
\begin{proof}
From (\ref{eq:b-sigma-bound}) we see that $|b'_U(x)|\leq \bar b$. 
Theorem 9, p.671 in \cite{aronson} and the subsequent remark state that the subtransition density $p_{\omega, U}(s,0,\cdot)$ is continuous in U.
From (\ref{eq:PDE-upper}) in Proposition \ref{prop:PDE} and from similar 
computations as carried out between (\ref{eq:sup-green}) and (\ref{eq:sup-green-1}), and applying dominated convergence, we see that $g_{U}(0,\cdot,\om)$ is continuous in  $U \smallsetminus \{0\}$. Consequently, by continuity of $b(\cdot,\om)$, see (\ref{eq:Lipschitz}),  and an application of (\ref{eq:green-upper}) and dominated convergence,  we see that $b'_U$ is continuous in $x \in U \smallsetminus \{0\}$.
\end{proof}
For $f \in C^{2}(\R{d})$, define
\begin{equation}
\label{eq:operator}
\mathcal{L}'f(x) \df \frac{1}{2} \Delta f(x) + b'_U(x) \nabla f(x)\,,
\end{equation}
and denote with, cf. \cite{bass} p.146, 
\begin{equation}
\label{eq:measure}
\text{$P'_{x,U}$ the unique solution to the martingale problem
for $\mathcal L'$ started at $x \in \R{d}$}.
\end{equation}
We write $E'_{x,U}$ for the corresponding expectation, and we denote with $p'_U(s,x,y)$, $x,y \in U$, $s>0$, 
the corresponding subtransition density (which can be defined by means of
Girsanov's theorem, see equation (4.1) in \cite{lyons}). 
Theorem 4.1 in \cite{lyons} states that estimate (\ref{eq:PDE-upper}) in Proposition \ref{prop:PDE} holds
for $p'_U$. With the same arguments as given in the proof of statement (\ref{eq:green-upper}) in Lemma \ref{cor:green}, we see that the Green function
\begin{equation}
  \label{eq:green-def}
  g'_U(x,y) \df \int_0^\infty p'_U(s,x,y)\textrm{d}s
\end{equation}
is well defined for $x,y \in U$, $x \neq y$, when $d \geq 2$, and for 
$x,y \in U$, when $d=1$.
The first step is
\begin{prop}
 \label{prop:exit}
 Let $U$ be a bounded $C^{\infty}$ domain containing 0. Then $X_{T_U}$ has same law under $P'_{0,U}$ and $P_{0}$ (see (\ref{eq:annealed}) for the
notation).
\end{prop}
\begin{proof}
We drop the subscript $U$ in $P'_{0,U}$ and $E'_{0,U}$.
By definition of the martingale problem, it holds for $f \in C^2(\R{d})$
that
\begin{equation*}
  E'_0[f(X_{t \wedge T_{U}})]-f(0)=E'_0\left[\int_{0}^{t \wedge T_{U}}
 \mathcal{L}'f(X_{s})\textrm{d}s\right]\,.
\end{equation*}
In particular, for $f \in C^{2}(\bar U)$, it follows from $E'_0[T_U]<\infty$ and from dominated convergence that
\begin{multline}
\label{eq:forward1}
  E'_0[f(X_{T_U})]=f(0)+E'_0[\int_{0}^{T_U}\mathcal{L}'f(X_{s})\textrm{d}s]\\
 =f(0)+\int_0^\infty E'_0[\mathcal{L}'f(X_{s}),\,s<T_U]\textrm{d}s
 =f(0)+\int_U g'_U(0,x)\mathcal{L}'f(x)\textrm{d}x\,.
\end{multline}
In the same way it follows that for $\om \in \Omega$,
\begin{equation}
  \label{eq:forward2}
  E_{0,\om}[f(X_{T_U})]=f(0)+\int_U g_U(0,x,\om)\mathcal{L}_\om f(x)\textrm{d}x\,. 
\end{equation}
Integrating (\ref{eq:forward2}) with respect to $\mathbb P$, the definition of $\mathcal L'$ (recall (\ref{eq:operator})) shows that
\begin{equation}
  \label{eq:forward3}
  E_0[f(X_{T_U})]=f(0)+\int_U \mathbb E [g_U(0,x,\om)]\mathcal{L}'f(x)\textrm{d}x\,.
\end{equation}
Combining (\ref{eq:forward1}) and (\ref{eq:forward3}), we obtain that for $f \in C^2(\bar U)$
\begin{equation}
  \label{eq:forward4}
 E_0[f(X_{T_U})]-E'_0[f(X_{T_U})]=\int_U (\mathbb E [g_U(0,x,\om)]-g'_U(0,x) )\mathcal{L}' f(x)\textrm{d}x\,. 
\end{equation}
Given $\phi \in C^\infty (\bar U)$, we will now find functions $u_n \in C^2(\bar U)$ such that 
\begin{equation}
\label{eq:smooth}
\lim_{n \to \infty}\mathcal L' u_n(x) =0 \text{ for a.e. } x \in U, \text{ and }
u_n=\phi  \text{ on the boundary } \partial U\,.
\end{equation}
Choose functions $b'_{U,n} \in  C^\infty (\bar U)$, $n \geq 1$, which converge 
boundedly a.e. in $U$ to $b'_U$. For $\phi \in  C^\infty (\bar U)$, consider the Dirichlet problem 
\begin{equation}
  \label{eq:dirichlet1}
  \frac{1}{2}\Delta u_n + b'_{U,n} \nabla u_n=0 \text{ in }U,\,u_n=\phi \text{ on }  \partial U\,.
\end{equation}
Following theorem 6.14 p.107 in \cite{gil-tru}, there is a unique solution
$u_n$ in $C^2(\bar U)$.  Fix $p>d$. The generalized problem
\begin{equation}
  \label{eq:dirichlet2}
  \mathcal L' u=0 \text{ in }U\,,\,u-\phi \in W_0 ^{1,p}(U)
\end{equation}
has a unique solution $u$ in the Sobolev space $W^{2,p}(U)$, see \cite{gil-tru} 
p.241.
Continuing our proof of (\ref{eq:smooth}), we will now show that 
\begin{equation}
  \label{eq:gradient-bounded}
  \sup_{n} \sup_{x \in U} |\nabla u_n (x)|<\infty\,.
\end{equation}
Define $w_n \df u_n-u$, $n \geq 1$, and obtain by means of the Sobolev inequality, see \cite{gil-tru} p.158, that
\begin{equation}
  \label{eq:sobolev}
  \sup_{x \in U}|\nabla u_n (x)| \leq  \sup_{x \in U}|\nabla w_n (x)| +  \sup_{x \in U}|\nabla u (x)|
\leq c(p,U) (\|w_n\|_{W^{2,p}(U)}+\|u\|_{W^{2,p}(U)})\,.
\end{equation}
$w_n$, $n \geq 1$, lies in the Sobolev space $W_0^{1,p}(U)$ and solves (see (\ref{eq:dirichlet1}) and (\ref{eq:dirichlet2}))
\begin{equation}
  \frac{1}{2}\Delta w_n + b'_{U,n} \nabla w_n=(b'_U-b'_{U,n})\nabla u \text{ in } U\,.
\end{equation}
Lemma 9.17 p.242 in \cite{gil-tru} and dominated convergence show that
\begin{equation}
  \label{eq:w_n}
  \|w_n\|_{W^{2,p}(U)}\leq c(p,U) \|(b'_U-b'_{U,n}) \nabla u\|_{L^p(U)} \underset{n \to \infty} \longrightarrow 0\,.
\end{equation}
Combining (\ref{eq:sobolev}), (\ref{eq:w_n}) and (\ref{eq:dirichlet2}) yields
(\ref{eq:gradient-bounded}). (\ref{eq:dirichlet1}) yields
\begin{equation}
  \mathcal L' u_n=(b'_U-b'_{U,n})\nabla u_n \text{ in }U,\,\,
  u_n=\phi  \text{ on } \partial U,
\end{equation}
which, together with (\ref{eq:gradient-bounded}), 
shows (\ref{eq:smooth}). Choosing $f=u_n$ in (\ref{eq:forward4}) and applying dominated convergence gives
\begin{equation}
  E_0[\phi(X_{T_U})]=E'_0[\phi(X_{T_U})] \text{ for all } \phi \in C^\infty (\bar U)\,.
\end{equation}
Since every function in $C^\infty(\partial U)$ is the restriction of a function in $C^\infty(\bar U)$,
see Lemma 6.37 p.137 in \cite {gil-tru}, the claim of the Proposition follows.
\end{proof}
We now introduce condition $(K)$, and show that it implies condition $(T)$. 
\begin{defn}
\label{def:K}
Let $l \in S^{d-1}$. We say that condition ($K)|l$ holds, if there is an
$\epsilon >0$, such that for all bounded domains $U$
containing 0 
\begin{equation}
\label{eq:K}
 \inf_{x \in U \smallsetminus \{0\}, \text{ dist} (x,\partial U)>5R} b'_U (x)\cdot l > \epsilon \,,
\end{equation}
with the convention $\inf \varnothing = +\infty$.
\end{defn}
\begin{prop}
\label{prop:K}
$(K)|l \Rightarrow (T)|l$\, (recall (\ref{eq:T})).
\end{prop}
\begin{proof}
The set of $l \in S^{d-1}$ for which (\ref{eq:K}) holds is open and  hence our claim will follow if for such an $l$ we show that
\begin{equation}
  \label{eq:T|l}
  \limsup_{L \to \infty}L^{-1}\log P_0[X_{T_{U_{l,b,L}}}\cdot l <0]<0\,.
\end{equation}
Denote with $\Pi_l(w) \df w-(w\cdot l) l$, $w \in \R{d}$, the projection on the orthogonal complement of $l$, and define 
\begin{equation}
  \label{eq:bounded-set}
  V_{l,b,L} \df \left\{x \in \R{d}: -bL < x \cdot l < L, |\Pi_l(x)| < L^{2} \right\}\,.
\end{equation}
In view of Proposition \ref{prop:exit}, we choose bounded  $C^\infty$ 
domains
$\tilde  V_{l,b,L}$ such that 
\begin{equation}
  \label{eq:set-inclusion}
   V_{l,b,L} \subset \left\{x \in \R{d}: -bL < x \cdot l < L, |\Pi_l(x)| < L^{2} +5R \right\} \subset \tilde  V_{l,b,L} \subset U_{l,b,L}\,.
\end{equation}
(When $d=1$, $\Pi_l(w) \equiv 0$, and we simply have that $U_{l,b,L}=V_{l,b,L}=\tilde V_{l,b,L}$.)
Recall (\ref{eq:measure}).
To prove (\ref{eq:T|l}), it will suffice to prove that
\begin{equation}
  \label{eq:exit-bounded-set}
  \limsup_{L \to \infty}L^{-1}\log  P'_{0,\tilde V_{l,b,L}}[X_{T_{V_{l,b,L}}}
  \cdot l <L]<0\,.
\end{equation}
Indeed, once this is proved, it follows from (\ref{eq:set-inclusion}) that
\begin{equation}
  \label{eq:exit-bounded-set2}
  \limsup_{L \to \infty}L^{-1}\log  P'_{0,\tilde V_{l,b,L}}[X_{T_{\tilde V_{l,b,L}}}\cdot l <L]<0\,.
\end{equation}
Hence, with Proposition \ref{prop:exit}, statement (\ref{eq:exit-bounded-set2})
holds with $ P'_{0,\tilde V_{l,b,L}}$ replaced by $P_0$, and, using (\ref{eq:set-inclusion}) once
more, (\ref{eq:T|l}) follows.\\
We now prove (\ref{eq:exit-bounded-set}). 
By (\ref{eq:set-inclusion}) and (\ref{eq:K}), we see that for $x \in V_{l,b,L}$,
\begin{equation}
  \label{eq:b'-bounds}
  b'_{\tilde V_{l,b,L}}(x)\cdot l \geq 
\begin{cases}
 \epsilon,  &\text{ if } -bL+5R<x \cdot l <L-5R \text{ and } x \neq 0,\\ 
 -\bar b,   &\text{ else }.
\end{cases}
\end{equation}
We thus consider the  process $X_t \cdot l$. 
We introduce the function $u(\cdot)$ on $\R{}$, which is
defined on $[-bL,L]$ through
\begin{equation}
u(r) \df
\begin{cases}
  \alpha_1 e^{\alpha_2 \epsilon (bL-5R)}(\alpha_3-e^{4\bar b (r-(-bL+5R))}),
  &\text{if } r \in [-bL,-bL+5R]\,,\\
  e^{-\alpha_2 \epsilon r}, &\text{if } r \in (-bL+5R,L-5R)\,,\\
  \alpha_4 e^{-\alpha_2 \epsilon (L-5R)}(\alpha_5-e^{4\bar b (r-(L-5R))}),
  &\text{if } r \in [L-5R,L]\,,
\end{cases}  
\end{equation}
and which is extended boundedly and in a $C^2$ fashion outside $[-bL,L]$, and such that $u$ is twice differentiable in the points $-bL$ and $L$.
The numbers $\alpha_i$, $1\leq i \leq 5$, are chosen positive and 
independent of $L$, via 
\begin{equation}
  \label{eq:coefficients}
  \alpha_5=1+e^{20\bar b R},\, \alpha_4=e^{-20\bar b R},\, \alpha_2 = \min (1, \frac{4 \bar b}{\epsilon}e^{-20\bar b R}),\,
 \alpha_1=\frac{\epsilon \alpha_2}{4 \bar b},\, \alpha_3=1+\frac{4 \bar b}{\epsilon \alpha_2}\,.
\end{equation}
Then, on $[-bL,L]$,  $u$ is positive, continuous and decreasing.
In addition, one has with the definition $j(r)=u'(r_+)-u'(r_-)$, 
\begin{equation}
\label{eq:derivative}
j(-bL+5R)=0, \text{ and } j(L-5R) \le 0\,.
\end{equation}
On $\R{d}$ we define the function $\tilde u(x)=u(x \cdot l)$, and 
for $\lambda$ real, we define on $\R{}_+ \times \R{}$ the function
$v_\lambda(t,r) \df e^{\lambda t} u(r)$, and on
$\R{}_+ \times \R{d}$ the function
 $\tilde v_\lambda(t,x) \df v_\lambda(t,x \cdot l)=e^{\lambda t} \tilde u(x)$. 
We will now find $\lambda_0$ positive such that
\begin{equation}
  \label{eq:supermartingale}
  v_{\lambda_0}(t \wedge T_{V_{l,b,L}}, X_{t \wedge T_{V_{l,b,L}}}\cdot l) \text{ is a positive supermartingale under }P'_{0,\tilde V_{l,b,L}}.
\end{equation}
Corollary 4.8 p.317 in \cite{kar-shr}, combined with remark 4.3 p.173 therein, shows the existence of a $d$-dimensional Brownian motion $W_t$ defined on 
$(C(\R{}_+,\R{d}), \mathcal F,P'_{0,\tilde V_{l,b,L}})$, such that
\begin{equation*}
 P'_{0,\tilde V_{l,b,L}}-\text{a.s.},\quad Y_t \df X_t \cdot l=
 W_t \cdot l+\int_0^t b'_{\tilde V_{l,b,L}}(X_s)\cdot l\,ds\,.
\end{equation*}
Writing $u$ as a linear combination of convex functions, we find from the generalised It\^o rule, see \cite{kar-shr} p.218, that 
 \begin{align}
  \label{eq:ito}
P'_{0,\tilde V_{l,b,L}}-\text{a.s.,}\quad u(Y_{t})\,=\,1+\int_0^{t}D^-u(Y_s)dY_s + \int_{-\infty}^{\infty}\Lambda_t(a)\mu(da),
\end{align}
where $D^-u$ is the left-hand derivative of $u$, $\Lambda(a)$ is the local time of $Y$ in $a$, and $\mu$ is the second derivative measure, i.e. 
$\mu([a,b))=D^-u(b)-D^-u(a)$, $a<b$ real. 
Notice that the first derivative  of $u$ exists and is continuous outside $L-5R$, and the second derivative of $u$ exists (in
particular) outside the Lebesgue zero set $A= \{-bL+5R,0,L-5R\}$.
Hence we find by definition of the second derivative measure, and  with the help of equation (7.3) p.218 in \cite{kar-shr} that
$P'_{0,\tilde V_{l,b,L}}$-a.s.,
\begin{equation}
\label{eq:sec-der}
\begin{aligned}
\int_{-\infty}^{\infty}\Lambda_t(a)\mu(da)=& \int_{-\infty}^{\infty}
\Lambda_t(a){\bf 1}_{A^c}(a) u''(a)\,da +\Lambda_t(L-5R)\,j(L-5R)\\
=&\tfrac{1}{2} \int_0^t  u''(Y_s){\bf 1}_{A^c}(Y_s)\,ds+\Lambda_t(L-5R)\,j(L-5R)\,.
\end{aligned}
\end{equation}
Another application of  equation (7.3) p.218 in \cite{kar-shr} shows that 
\begin{equation}
\label{eq:N}
P'_{0,\tilde V_{l,b,L}}-\text{a.s.,}\quad  \int_0^t {\bf 1}_{A}(Y_s)\,ds=2\int_{-\infty}^\infty {\bf 1}_{A}(a)\Lambda_t(a)\,da=0\,.
\end{equation}
As a result, we find that $P'_{0,\tilde V_{l,b,L}}$-a.s.,
\begin{equation}
\label{eq:first-der}
  \int_0^t D^-u(Y_s){\bf 1}_{A}(Y_s)dY_s =0\,.
\end{equation}
Combining (\ref{eq:sec-der}) and (\ref{eq:first-der}), and by definition of the operator $\mathcal L'$, see (\ref{eq:operator}), 
we can now rewrite (\ref{eq:ito}) as the
$P'_{0,\tilde V_{l,b,L}}$-a.s. equalities
\begin{equation*}
\begin{aligned}
u(Y_{t})=&1+\int_0^{t} u'(Y_s){\bf 1}_{A^c}(Y_s)dY_s 
+\tfrac{1}{2}\int_0^{t} u''(Y_s) {\bf 1}_{A^c}(Y_s)ds + 
\Lambda_t(L-5R)j(L-5R)\\
=&1+\int_0^{t}\mathcal L'\tilde u(X_s){\bf 1}_{A^c}(X_s \cdot l)\,ds + \Lambda_t(L-5R)\,j(L-5R)
+M_t\,,
\end{aligned}
\end{equation*}
where $M_t$ is a continuous martingale.
In particular, 
$\tilde u(X_t) (=u(Y_t))$ is a continuous semimartingale, and 
applying It\^o's rule to the product $e^{\lambda t} \cdot \tilde u(X_t)=\tilde v_\lambda(t,X_t)$, and using (\ref{eq:N}) once again, we obtain that,
$P'_{0,\tilde V_{l,b,L}}$-a.s.,
\begin{equation}
\label{eq:ito-2}
\begin{aligned}
&\tilde v_{\lambda}(t,X_{t })
=\,1+\int_0^{t }\lambda e^{\lambda s}
  \tilde u(X_s) \,ds + \int_0^{t } e^{\lambda s}\, d\tilde u(X_s)\\
&=\,1+ \int_0^{t } \Big(\tfrac{\partial}{\partial s}+
   \mathcal L'\Big)\tilde v_\lambda (s,X_s){\bf 1}_{A^c}(X_s \cdot l)\,ds
 +j(L-5R)\int_0^{t } e^{\lambda s}d\Lambda^{L-5R}_s +N_{t }\,,
\end{aligned}
\end{equation}
where $N_{t }$ is a continuous martingale.
We find through direct computation that for $ x \in V_{l,b,L}$, and a suitable $\psi (x) \geq 0$, using the notation $I_1=(-bL,-bL+5R)$, $I_2=(-bL+5R,L-5R)$, $I_3=(L-5R,L)$,
\begin{equation*}
 \big[(\tfrac{\partial}{\partial s}+ \mathcal L')\tilde v_\lambda \big](s,x)\leq \psi (x)e^{\lambda s} \cdot
\begin{cases}
\lambda(e^{20\bar b R}\alpha_3-1)-4\bar b( 2 \bar b + b'_{\tilde V_{l,b,L}} (x) \cdot l)\,,&\text{ if } x \cdot l \in I_1\,,\\ 
\lambda+\alpha_2\epsilon(\frac{1}{2}\alpha_2\epsilon-b'_{\tilde V_{l,b,L}} (x) \cdot l)\,,&\text{ if }x \cdot l \in I_2\,,\\
\lambda(\alpha_5-1)-4\bar b( 2 \bar b + b'_{\tilde V_{l,b,L}} (x) \cdot l)\,,&\text{ if }x \cdot l \in I_3\,.
\end{cases}
\end{equation*}
Hence, by (\ref{eq:b'-bounds}) and (\ref{eq:coefficients}), we can find 
$\lambda_0>0$ small such that for $x \in V_{l,b,L}$, $x \cdot l \notin A$, the right-hand side of the last expression is negative. 
Since $j(L-5R) \le 0$, see (\ref{eq:derivative}), we obtain from (\ref{eq:ito-2}) applied to the 
finite stopping time $t \wedge T_{V_{l,b,L}}$ that (\ref{eq:supermartingale}) holds.\\
We now derive the claim of the proposition from (\ref{eq:supermartingale}).
When $d \geq 2$, the probability to exit $V_{l,b,L}$ neither from the ``right'' nor from the ``left''
can be bounded as follows:
\begin{equation}
\label{eq:bound}
\begin{aligned}
  & P'_{0,\tilde V_{l,b,L}}[-bL<X_{T_{V_{l,b,L}}} \cdot l < L\,]\leq \\
 & P'_{0,\tilde V_{l,b,L}}[-bL<X_{T_{V_{l,b,L}}} \cdot l < L,\,T_{V_{l,b,L}} > 
\tfrac{2\alpha_2 \epsilon}{\lambda_0}L\,]
+ P'_{0,\tilde V_{l,b,L}}[\,\sup |X_t|\geq L^2:t \leq \tfrac{2\alpha_2 \epsilon}{\lambda_0}L\,]\,.
\end{aligned}
\end{equation}
By Chebychev's inequality and Fatou's lemma, we find that the first term on the 
right-hand side is smaller than
\begin{equation}
\label{eq:bound-0}
\begin{aligned}
&\frac{1}{v_{\lambda_0}(\tfrac{2\alpha_2 \epsilon}{\lambda_0}L,L)}E'_{0,\tilde V_{l,b,L}}
[v_{\lambda_0}(T_{V_{l,b,L}}, X_{T_{V_{l,b,L}}}\cdot l)]\\
\le &c(\epsilon)e^{-\alpha_2 \epsilon L}   \liminf_{t \to \infty}E'_{0,\tilde V_{l,b,L}}
[v_{\lambda_0}(t \wedge T_{V_{l,b,L}}, X_{t \wedge T_{V_{l,b,L}}}\cdot l)]\\
\le & c(\epsilon)e^{-\alpha_2 \epsilon L}\,   v_{\lambda_0}(0,0)=c(\epsilon)e^{-\alpha_2 \epsilon L},
\end{aligned}
\end{equation}
where, in the last inequality, we used (\ref{eq:supermartingale}). 
Applying (\ref{eq:max}) in Lemma \ref{lemma:bernstein} to the second term in the right-hand side of (\ref{eq:bound}), we obtain, together with (\ref{eq:bound-0}), that
\begin{equation}
\label{eq:bound-1}
\limsup_{L \to \infty}L^{-1} \log P'_{0,\tilde V_{l,b,L}}[-bL<X_{T_{V_{l,b,L}}} \cdot l < L\,]<0.
\end{equation}
When $d \geq 1$, we bound the probability to exit $V_{l,b,L}$ from the left
by a similar argument as in (\ref{eq:bound-0}), and find that
\begin{equation}
  \label{eq:bound-2}
   P'_{0,\tilde V_{l,b,L}}[X_{T_{V_{l,b,L}}} \cdot l=-bL] \leq \frac{v_{\lambda_0}(0,0)}{v_{\lambda_0}(0,-bL)}
\leq e^{-c(\epsilon)L}\,.
\end{equation}
(\ref{eq:bound-2}), together with  (\ref{eq:bound-1}), when $d \geq 2$, show 
(\ref{eq:exit-bounded-set}), which implies condition $(T)|l$.
\end{proof}
Let us now turn to the
\begin{proof}[{\bf Proof of Theorem \ref{thm:examples}}]
It suffices to verify condition $(K)|l$, which implies condition $(T)|l$, see Proposition
\ref{prop:K}. Let $U$ be a bounded domain containing 0, and assume that there is
\begin{equation}
  \label{eq:dist}
  x \in U \smallsetminus \{0\} \text{ such that }\text{dist}(x,\partial U)>5R\,.
\end{equation}
(otherwise $(K)|l$ automatically holds).
With $x$ as above, $\delta >0$, for $f$ a non-negative bounded measurable function on $U$, we write
\begin{align*}
f_\delta(\cdot) \df f(\cdot) \mathbf{1}_{B_{\delta}(x)}(\cdot)\,,\text{ and }
b^{\pm}_\delta(\cdot,\om) \df (b(\cdot,\om)\cdot l)_{\pm} \mathbf{1}_{B_{\delta}(x)}(\cdot)\,.  
\end{align*}
Lemma \ref{lemma:b'} shows that
\begin{equation}
\label{eq:b'-integral}
 b'_U(x)\cdot l=\lim_{\delta \to 0}\frac{1}{|B_\delta|}\int_{B_{\delta}(x)}b'_U(y)\cdot l~\textrm{d}y\,. 
\end{equation}
If we choose $\delta < |x|/2$, it follows from (\ref{eq:green-lower})
and from (\ref{eq:green-upper}) in Corollary \ref{cor:green} that
\begin{equation}
\label{eq:inf-sup-green}
  0<\inf_{y \in B_{\delta}(x)}\mathbb{E}[g_{U}(0,y,\om)]\leq \sup_{y \in B_{\delta}(x)}\mathbb{E}[g_{U}(0,y,\om)]<\infty\,,
\end{equation}
and  we obtain by the definition of $b'_U$, see (\ref{eq:b'}), that
\begin{equation}
  \label{eq:b'-bound}
  \begin{aligned}
 \frac{1}{|B_\delta|}\underset{B_{\delta}(x)}{\int}b'_U(y)\cdot l~\textrm{d}y
 \geq  \frac{\mathbb{E}\left[\int g_{U}(0,y,\om)b_{\delta}^{+}(y,\om)
\textrm{d}y \right]}{|B_\delta|\,\sup_{y \in B_\delta(x)}\mathbb{E}[g_{U}(0,y,\om)]}-
\frac{\mathbb{E}\left[\int g_{U}(0,y,\om)b_{\delta}^{-}(y,\om)
\textrm{d}y \right]}{|B_\delta|\,\inf_{y \in B_\delta(x)}\mathbb{E}[g_{U}(0,y,\om)]}\,.
\end{aligned}
\end{equation}
Denote with $R_k$ and $D_k$, $k \geq 1$, the successive returns of $X_\cdot$ to $B_{2R}(x)$ and departures from 
$B_{4R}(x)$ defined similarly as in (\ref{eq:excursion-1}) and (\ref{eq:excursion-2}), with $B_1(x)$ and $B_2(x)$ replaced by $B_{2R}(x)$ and $B_{4R}(x)$ 
respectively. 
For $y$ in $U$, define the associated operators:
\begin{equation*}
 Rf(y) \df E_{y,\om}\left[f(X_{R_{1}}), R_{1} < T_{U} \right],\,
 Qf(y) \df E_{y,\om}\left[f(X_{D_{1}})\right],\,
 Tf(y) \df  E_{y,\om}[\int_{0}^{D_{1}}f(X_{s})\textrm{d}s].
\end{equation*}
If $\delta \leq 2R$, successive applications of the strong Markov property
show that
\begin{equation}
\label{eq:Markov}
\int_{U}g_{U}(0,y,\om)f_{\delta}(y)\textrm{d}y 
= E_{0,\om}[\int_{0}^{T_{U}}f_{\delta}(X_s)\textrm{d}s]
= R(Id - QR)^{-1}Tf_{\delta}(0).
\end{equation}
In view of (\ref{eq:b'-bound}), it will be crucial to bound the above quantity 
from below and from above. In a first step, we derive bounds
on the operators $R$ and $QR$. For $y \in U$, we have 
\begin{equation}
\label{eq:R-bound}
\inf_{z \in \partial{B}_{2R}(x)}f(z)P_{y,\om}[R_{1}<T_{U}]\,\,\leq Rf(y)\, \,
\leq \sup_{z \in \partial{B}_{2R}(x)}f(z)P_{y,\om}[R_{1}<T_{U}],
\end{equation}
and hence,
\begin{equation}
\label{eq:QR-bound}
\begin{aligned}
&\sup_{y \in \partial B_{2R}(x)}QRf(y) \leq \sup_{z \in \partial{B}_{4R}(x)}
  P_{z,\om}[R_{1}<T_{U}] \sup_{z \in \partial{B}_{2R}(x)}f(z),\\
&\inf_{y \in \partial B_{2R}(x)}QRf(y) \geq \inf_{z \in \partial{B}_{4R}(x)}
  P_{z,\om}[R_{1}<T_{U}] \inf_{z \in \partial{B}_{2R}(x)}f(z)\,.
\end{aligned}
\end{equation}
We first derive a lower bound for (\ref{eq:Markov}), see (\ref{eq:estimate-2})
below. Repeated applications of (\ref{eq:R-bound}) and (\ref{eq:QR-bound}) yield
\begin{equation}
  \label{eq:estimate-1}
\begin{aligned}
& R(Id - QR)^{-1}Tf_\delta(0)\\
\geq & \, P_{0,\om}[R_{1}<T_{U}]\sum_{j \geq 0}\left(\inf_{z \in \partial
          {B}_{4R}(x)}P_{z,\om}[R_{1}<T_{U}]\right)^{j}\inf_{z \in 
          \partial{B}_{2R}(x)}Tf_\delta(z)\\
\geq & \, \frac{ P_{0,\om}[R_{1}<T_{U}]}{\sup_{z \in \partial B_{4R}(x)}P_{z,\om}
          [R_{1}>T_{U}]}\inf_{z \in B_\delta (x)} f_\delta(z)
        \inf_{z \in \partial{B}_{2R}(x)}T \mathbf{1}_{B_{\delta}(x)}(z)\,.
\end{aligned}  
\end{equation}
If $\delta <2R$, we find by means of (\ref{eq:green-lower}) in Corollary \ref{cor:green} that 
\begin{equation}
\label{eq:green-1}
\inf_{z \in \partial{B}_{2R}(x)}T \mathbf{1}_{B_{\delta}(x)}(z)
 \geq  \int_{B_{\delta}(x)} ~\inf_{z \in \partial{B}_{2R}(x)}g_{B_{4R}(x)}(z,y,\om)\textrm{d}y 
\geq c\, |B_\delta|\,.
\end{equation}
Combining (\ref{eq:estimate-1}) and (\ref{eq:green-1}), and using (\ref{eq:Markov}), we see that 
\begin{equation}
  \label{eq:estimate-2}
\int_{U}g_{U}(0,y,\om)f_{\delta}(y)\textrm{d}y  
 \geq c\,|B_\delta|\,\frac{ P_{0,\om}[R_{1}<T_{U}]}{\sup_{z \in \partial  B_{4R}(x)}P_{z,\om}[R_{1}>T_{U}]}\inf_{z \in B_\delta (x)} f_\delta(z)\,.
\end{equation}
We will now derive an upper bound on (\ref{eq:Markov}), see (\ref{eq:estimate-3}). If $\delta < R$, we find by another use of Corollary \ref{cor:green} that 
\begin{equation}
 \sup_{z \in \partial{B}_{2R}(x), y \in B_\delta (x)}g_{B_{4R}(x)}(z,y,\om) \leq c\,.
\end{equation}
Proceeding in a similar fashion as in (\ref{eq:estimate-1})-(\ref{eq:estimate-2}), we obtain the upper bound
\begin{equation}
\label{eq:estimate-3}
\int_{U}g_{U}(0,y,\om)f_{\delta}(y)\textrm{d}y  
\leq c\,|B_\delta |\,\frac{ P_{0,\om}[R_{1}<T_{U}]}{\inf_{z \in \partial B_{4R}(x)}P_{z,\om}[R_{1}>T_{U}]}\,\sup_{z \in B_{\delta}(x)}f_\delta(z)\,.
\end{equation}
We will now give a lower bound for the first term in the last line of
(\ref{eq:b'-bound}).
Applying (\ref{eq:Markov}) with $f_\delta=b^+_\delta$  
and using (\ref{eq:estimate-2}), we see that 
\begin{equation}  
  \label{eq:estimate-4} 
\mathbb{E}\left[\int g_{U}(0,y,\om)b_{\delta}^{+}(y,\om)
\textrm{d}y \right]
\geq
 c\,|B_\delta |\,
\mathbb{E}\left[\frac{ P_{0,\om}[R_{1}<T_{U}]}{\sup_{z \in \partial B_{4R}(x)}P_{z,\om}[R_{1}>T_{U}]}\,\inf_{z \in B_\delta (x)} b_{\delta}^{+}(z,\om)\right]\,. 
\end{equation} 
Observe that $P_{0,\om}[R_{1}<T_U]=1$ if $0 \in B_{2R}(x)$. Hence 
$\frac{ P_{0,\om}[R_{1}<T_{U}]}{\inf_{z \in \partial B_{4R}(x)}P_{z,\om}[R_{1}
>T_{U}]}$ is $\mathcal{H}_{B_{2R}^{c}(x)}$-measurable. Since 
$\inf_{z \in B_{\delta}(x)}b_{\delta}^{+}(z,\om)$ is 
$\mathcal{H}_{B_{\delta}(x)}$-measurable, it follows for $\delta < R$ and from finite range dependence, see (\ref{eq:R-separation}), that these two random variables are $\mathbb{P}$-independent,  and hence (\ref{eq:estimate-4}) equals
\begin{equation}
  \label{eq:estimate-5}
  c\,|B_\delta |\, 
\mathbb{E}\left[\frac{ P_{0,\om}[R_{1}<T_{U}]}{\sup_{z \in \partial B_{4R}(x)}P_{z,\om}[R_{1}>T_{U}]}\right] \mathbb{E}[\inf_{z \in B_\delta (x)} b_{\delta}^{+}(z,\om)]\,. 
\end{equation}
The application of Harnack's inequality (see \cite{gil-tru} p.199) to the 
$\mathcal L_\om$-harmonic function $P_{\cdot,\om}[R_{1}>T_{U}]$ on $B_{5R}(x) \smallsetminus \bar{B}_{2R}(x)$ shows that 
\begin{equation*}
  \sup_{z \in \partial B_{4R}(x)}P_{z,\om}[R_{1}>T_{U}] \leq c\, \inf_{z \in \partial B_{4R}(x)}P_{z,\om}[R_{1}>T_{U}]\,.
\end{equation*}
Together with an application of (\ref{eq:estimate-3}) and (\ref{eq:Markov}) with $f_\delta=\mathbf 1_{B_\delta(x)}$, we obtain that (\ref{eq:estimate-5}) is 
bigger than
\begin{multline}
  \label{eq:estimate-6}
  c\,\mathbb{E}\left[\int_{B_{\delta}(x)}g_{U}(0,y,\om)
\textrm{d}y\right]\,\mathbb{E}[\inf_{z \in B_{\delta}(x)}b_{\delta}^{+}(z,\om)]\\
\geq c\,|B_\delta |\,\mathbb{E}[\inf_{y \in B_\delta (x)}
g_{U}(0,y,\om)]\,\mathbb{E}[\inf_{z \in B_{\delta}(x)}b_{\delta}^{+}(z,\om)]\,.
\end{multline}
Finally, using (\ref{eq:estimate-4})-(\ref{eq:estimate-6}), we find that the first term in the right-hand side of (\ref{eq:b'-bound}) is bigger than
\begin{equation}
  \label{eq:estimate-7}
 c_1\,\frac{\mathbb{E}[\inf_{y \in B_{\delta}(x)}g_{U}(0,y,\om)]}
{\mathbb{E}[\sup_{y \in B_{\delta}(x)}g_{U}(0,y,\om)]}
~\mathbb{E}[\inf_{z \in B_{\delta}(x)}b_{\delta}^{+}(z,\om)]\,.
\end{equation}
By similar computations as carried out between (\ref{eq:estimate-4}) and (\ref{eq:estimate-7}), we find as an upper bound for the second term in the right-hand side of (\ref{eq:b'-bound})
\begin{equation}
  \label{eq:estimate-8}
 c_2\, \frac{\mathbb{E}[\sup_{y \in B_{\delta}(x)}g_{U}(0,y,\om)]}
  {\mathbb{E}[\inf_{y \in B_{\delta}(x)}g_{U}(0,y,\om)]}
  ~\mathbb{E}[\sup_{z \in B_{\delta}(x)}b_{\delta}^{-}(z,\om)]\,.
\end{equation}
The continuity of $g_U(0,\cdot,\om)$ and of $b_\delta ^+ (\cdot, \om)$ in 
$B_\delta (x)$, see lemma \ref{lemma:b'} and (\ref{eq:Lipschitz}), together with  dominated convergence, and the translation invariance of the measure $\mathbb P$, show that
\begin{equation}
  \label{eq:estimate-9}
  \lim_{\delta \to 0}\,c_1\,\frac{\mathbb{E}[\inf_{y \in B_{\delta}(x)}g_{U}(0,y,\om)]}
{\mathbb{E}[\sup_{y \in B_{\delta}(x)}g_{U}(0,y,\om)]}
~\mathbb{E}[\inf_{z \in B_{\delta}(x)}b_{\delta}^{+}(z,\om)]
=\,c_1\,\mathbb{E}[(b(0,\om)\cdot l)_+]\,,
\end{equation}
and a similar identity for the term in (\ref{eq:estimate-8}). 
Inserting (\ref{eq:estimate-7})-(\ref{eq:estimate-9}) in (\ref{eq:b'-bound}),
and using (\ref{eq:b'-integral}), we finally obtain
\begin{equation}
\label{eq:estimate-10}
 b'_U(x)\cdot l \geq ~c_1~\mathbb{E}[(b(0,\om)\cdot l)_{+}- \tfrac{c_2}{c_1}\,(b(0,\om)\cdot l)_{-}]\,.
\end{equation}
Hence, if (\ref{eq:criterion}) holds with $c_e \df \tfrac{c_2}{c_1}$,
we see that there is an $\epsilon >0$ such that for all $x$ as in (\ref{eq:dist})
\begin{equation}
\label{eq:estimate-11}
   b'_U(x)\cdot l > \epsilon\,.
\end{equation}
We conclude that condition $(K)|l$ holds, see (\ref{eq:K}). By means of Proposition \ref{prop:K}, condition $(T)|l$ holds, and 
Theorem \ref{thm:examples} is proved.
\end{proof}
\begin{rem} \rm
\label{rem:examples}
With the help of Theorem \ref{thm:examples}, it is easy
to obtain concrete examples of diffusions fulfilling condition $(T)$. 
For instance, when $(b(0,\om)\cdot l)_- =0$, we find:
\begin{equation}
\label{eq:non-nestling}
\begin{aligned}
&\text{Condition $(T)$ holds when $d \geq 1$ and there is $l \in S^{d-1}$ and $\delta >0$,} \\
&\text{such that $b(0,\om)\cdot l \geq 0$ for all $\om \in \Omega$, and  
 $p_\delta = \mathbb P[b(0,\om)\cdot l \geq \delta]>0$}\,.
\end{aligned}
\end{equation}
If there is $\delta >0$ such that $p_\delta=1$, this is in the spirit of the 
{\it non-nestling} case, which is in fact already covered by Proposition
\ref{prop:non-nestling}, and else, of the {\it marginal nestling} case
in the discrete setting,  see Sznitman \cite{szn00}.

Of course, Theorem \ref{thm:examples} also comprises more involved examples of condition $(T)$ where \\
$b(0,\om)\cdot l$ takes both positive and negative values for every $l \in S^{d-1}$.
Hence, when $d \geq 2$, Theorem \ref{thm:examples} provides examples of ballistic diffusions in random environment beyond previous knowledge. They  correspond to the {\it plain nestling} case in \cite{szn00}. 
\end{rem}

\section{Appendix}
\label{sec:appendix}
\small
\subsection{Bernstein's Inequality}
Recall the convention of the constants stated at the end of the Introduction. The following Lemma follows in essence from Bernstein's inequality (see \cite {rev-yor} page 153-154).
\begin{lemma}
\label{lemma:bernstein}
On $\R{d}$ we consider measurable functions $a$, $b$, with values in the space of symmetric matrices
and in $\R{d}$ respectively, that satisfy for suitable $\nu \ge 1$, and $\bar a >0$,  $\bar b >0$, 
\begin{equation}
  \label{eq:1}
  \tfrac{1}{\nu}|y|^{2}\leq \sum_{i,j} a_{ij}(x)y_i y_j \leq \nu |y|^{2},\,\,      \left|a(x) \right| \leq \bar a,\,\,
  \left|b(x) \right| \leq \bar b,\,\,x,y \in \R{d}\,.
\end{equation}
We denote with $\mathcal L$ the operator attached to $a$ and $b$, similarly as in (\ref{eq:diff-operator}), and we assume that $P_x$ solves the martingale problem for $\mathcal L$ started at $x$ in $\R{d}$. We denote with $E_x$ the corresponding expectation. Write $(X_t)_{t \ge 0}$ for the canonical process on
$C([0,\infty),\R{d})$, and let $Z_t = \sup_{s \leq t}|X_s-X_0|$.
Then, for every $\alpha >0$, 
 there are two constants $c(\alpha)>0$ and  $\tilde c(\alpha)>0$, such that 
for large $L$,
  \begin{equation}
    \label{eq:max}
    \sup_{x}P_x \big[Z_{\alpha L}\geq L^2\big]
    \leq \tilde c e^{-cL^3}\,.
  \end{equation}
Further, for $\gamma \in (0,1]$ and for all $\alpha>0$, there exists a constant $\delta(\alpha)>0$ such that
  \begin{equation}
    \label{eq:Z}
    \sup_{x} E_{x}\big[e^{\delta Z_1^\gamma}\big]\leq 1+\alpha\,.
  \end{equation}
\end{lemma}
\begin{proof}
We obtain from the martingale problem that $M_t=X_t-X_0-\int_0^t b\,(X_s)\,ds$ is a martingale. We compute the bracket  $\langle M^i \rangle_t$ of the $i$-th component $M^i_t$ of $M_t$, $1 \le i \le d$, and find  
$\langle M^i \rangle_t=\int_0^t a_{ii}(X_s)\,ds$. (\ref{eq:1}) yields 
$\langle M^i \rangle_t \le \nu t$, and with the help of  Bernstein's inequality (see \cite {rev-yor} page 153-154) and a further application of (\ref{eq:1}),
it follows immediately that for large $L$,
\begin{equation*}
  P_{x}\big[Z_{\alpha L}\geq L^2\big]
 \leq P_{x}\big[\sup_{s\leq \alpha L}|M_s|\geq (L^2-\alpha\bar b L) \big] 
 \leq 2d e^{-\frac{L^3}{4\nu \alpha d}}\,,
\end{equation*}
which proves (\ref{eq:max}).
Since $Z_1\leq \sup_{s\leq 1}|M_s|+\bar b$, we obtain for $0<\delta<1$ that
\begin{align*}
  &E_{x}\big[e^{\delta Z_1^\gamma}\big]\leq e^{\delta \bar{b}^\gamma}\, E_{x}\big[\exp\{\delta (\sup_{s\leq 1}|M_s|)^\gamma\}\big]\\
  =&e^{\delta\bar b^\gamma}\Big(1+ \delta\int^\infty_{0} \mathrm{d}v\; e^{\delta v}\underbrace{P_{x}\big[(\sup_{s\leq 1}|M_s|)^\gamma \geq v\big]}_{\leq 2d \exp\{-v^{\frac{2}{\gamma}}/(2d\nu)\}}\Big) 
  \leq e^{\delta\bar b^\gamma}\big(1+\delta\, c \big)\,,
\end{align*}
which proves (\ref{eq:Z}).
\end{proof}

\subsection{Bounds on the Green function}
The bounds on the transition density contained in the next Proposition will be crucial to derive bounds on the Green function.
\begin{prop}
  \label{prop:PDE}
Let $\mathcal L_\om$ be as in (\ref{eq:diff-operator}), and let assumptions
(\ref{eq:b-sigma-bound})-(\ref{eq:elliptic}) be in force. Then the linear parabolic equation of second order $\frac{\partial u}{\partial t}=\mathcal L_\om u$ has a  a unique fundamental solution $p_\om(t,x,y)$, and there are positive constants $\alpha$, $\beta$, $a$ and $\tilde \alpha$ such that for $t \leq 1$
\begin{equation}
    \label{eq:PDE-upper}
    |p_\om(t,x,y)| \leq \frac{\alpha}{t^{d/2}} \exp\big\{-\tfrac{\beta |x-y|^2}{t}\big\}\,,
  \end{equation}
and such that for  $|x-y|^2< a t$ and $t\in (0, 1]$
\begin{equation}
    \label{eq:PDE-lower}
    p_\om(t,x,y)\geq  \frac{\tilde \alpha}{t^{d/2}}\,.
  \end{equation}
\end{prop}
For the proof we refer the reader to \cite{illin}. The statements (4.16) and (4.75) therein correspond to (\ref{eq:PDE-upper}) and (\ref{eq:PDE-lower}).
Recall the convention on the constants stated at the end of the Introduction.
We obtain the following Corollary:
\begin{cor}
 \label{cor:green}
Assume (\ref{eq:b-sigma-bound}) and (\ref{eq:Lipschitz}), 
 and let $U$ be a bounded domain. There is a positive constant $m(r,U)$ such that for all $\om \in \Omega$, and for all $y, z \in U$ with dist$(y, \partial U)>r$, dist$(z, \partial U)>r$,
\begin{equation}
 \label{eq:green-lower}
  g_{U}(y,z,\om) \geq m\,.
\end{equation}
For $y \neq z$, define 
\begin{equation}
\label{eq:function-h}
h_y(z)=  
  \begin{cases}
   |y-z|^{2-d} \,, &d \geq 3\,,\\
   \log \frac{\text{diam}(U)}{|y-z|}\,,&d=2\,.
  \end{cases}
\end{equation}
There are positive constants $\alpha, c(U)$ such that for $y,z \in U$, and all
$\om \in \Omega$,
\begin{equation}
  \label{eq:green-upper}
  g_U(y,z,\om) \leq 
\begin{cases}
\alpha h_y(z)+c, & \text{if $d \geq 2$ and $y \neq z$},\\
c, &\text{if $d=1$}\,.
\end{cases}
\end{equation}
\end{cor}
\begin{proof}
Let $x \in U$ with dist$(x, \partial U)>r$. Choose $t_{0} \in (0,1]$ such that $\sqrt{at_0}\leq \frac{r}{2}$ and  for all $t\leq t_0$, $\frac{\tilde \alpha}{t^{d/2}}\geq \frac{2\alpha}{t_0^{d/2}} \exp\{-\tfrac{\beta r^2}{4t_0}\}$ holds, and such that in addition the function $t\mapsto \frac{\alpha}{t^{d/2}} \exp\{-\tfrac{\beta r^2}{4t}\}$ is monotone increasing on $\{t:t\leq t_0\}$.
Let $\rho =\min(\frac{r}{2},\sqrt{at_{0}/8})$ and $z_0 \in B_{\rho}(x)$. Hence 
$|X_{T_U}-z_0|>\tfrac{r}{2}$, and on the event $\{T_U<t\leq t_0\}$, the inequality $p_\omega(t-T_U, X_{T_U}, z)\leq \frac{\alpha}{t^{d/2}} \exp\big\{-\tfrac{\beta r^2}{4t}\big\}$ follows from (\ref{eq:PDE-upper}) and from the monotonicity mentioned above. Choose further $y_0 \in B_\rho(x)$, then $|y_0-z_0|<\sqrt{at_{0}/2}$, and hence, for $t \in (t_0/2,t_0)$, $|y_0-z_0|<\sqrt{at}$ holds.
By Duhamel's formula, see  \cite{stroock} page 331, and by (\ref{eq:PDE-lower}), the subtransition density $p_{\om,U}(t,y,z)$ satisfies for $y_0,z_0 \in B_{\rho}(x)$ and $t \in (t_0/2,t_0)$ 
\begin{equation}
\label{eq:inf-green}
\begin{aligned}
  p_{\om,U}(t,y_0,z_0)&=p_\om(t,y_0,z_0)-E_{y_0,\om}\big[T_U<t, p_\om(t-T_U, X_{T_U}, z_0)\big]\\
  &\geq \tfrac{\alpha}{t_0^{d/2}}\exp\{-\tfrac{\beta r^2}{4t_0}\}>0\,.
\end{aligned}
\end{equation}
We will now prove (\ref{eq:green-lower}).
Since $U$ is a bounded domain, it follows from a standard chaining argument using (\ref{eq:inf-green}) that there is a finite integer $K(U)>0$ such that for all $y,z \in U$ as above (\ref{eq:green-lower}), for all $t \in (Kt_0/2,Kt_0)$ and for all $\om \in \Omega$, 
\begin{equation}
  p_{\om,U}(t,y,z) \geq c(r,K)>0\,.
\end{equation}
Since 
\begin{equation*}
  g_{U}(y,z,\om) \geq  \int_{K\frac{t_0}{2}}^{K t_0}p_{\om,U}(t,y,z)\textrm{d}t\,,
\end{equation*}
the claim (\ref{eq:green-lower}) follows.
To prove the upper bound (\ref{eq:green-upper}), we write
\begin{align}
\label{eq:sup-green}
  g_{U}(y,z,\om)=\int_{0}^{\infty}p_{\om,U}(t, y  ,z)\textrm{d}t
   \leq \int_{0}^{1}p_{\om}(t, y  ,z)\textrm{d}t + \sum_{k=2}^{\infty}\int_\frac{k}{2}^{\frac{k+1}{2}}p_{\om,U}(t,y,z)\textrm{d}t\,.
\end{align}
With the help of (\ref{eq:PDE-upper}), we find positive constants $\alpha,\,c$
 such that 
\begin{equation}
  \label{eq:h}
  \int_{0}^{1}p_{\om}(t, y  ,z)\textrm{d}t \leq 
  \begin{cases}
  \alpha h_y(z)+c\,,& \text{ if $d \geq 2$, $y \neq z$},\\
  c\,, &\text{ if $d=1$}\,.
  \end{cases}
\end{equation}
We obtain by a repeated use of the Chapman-Kolmogorov equation and by (\ref{eq:PDE-upper}), that for $k \geq 2$,
\begin{multline*}
  \int_\frac{k}{2}^{\frac{k+1}{2}}p_{\om,U}(t,y,z)\textrm{d}t
  \leq  \int_U \textrm{d}v ~p_{\om,U}(1/2,y,v) ~\sup_{v \in U}\int_\frac{  k-1}{2}^{\frac{k}{2}}p_{\om,U}(t,v,z)\textrm{d}t\\
  \stackrel{induction}{\leq}  \left(\sup_{v \in U}P_{v,\om}[T_{U}>\frac{1}{2}]  \right)^{k-1}
  \sup_{v \in U} \int_{\frac{1}{2}}^{1}p_{\om,U}(t,v,z)\textrm{d}t
  \leq c \left(\sup_{v \in U}P_{v,\om}[T_{U}>\frac{1}{2}]\right)^{k-1}.
\end{multline*}
Hence, with the help of the Support Theorem of Stroock-Varadhan, see \cite{bass} p.25, or from a chaining argument using (\ref{eq:PDE-lower}), the sum on the right-hand side of (\ref{eq:sup-green}) will be smaller than
\begin{equation}
\label{eq:sup-green-1}
   \frac{c}{\inf_{v \in U}P_{v,\om}
  [T_{U} \leq\frac{1}{2}]} \leq  c(U) < \infty\,.
\end{equation}
Combining (\ref{eq:sup-green}), (\ref{eq:h}) and (\ref{eq:sup-green-1}) shows  (\ref{eq:green-upper}). 
\end{proof}


\begin{thebibliography}{30}

\bibitem{aronson}Aronson, D.G.: ``Non-negative solutions of linear parabolic equations'', {\it Ann. Scuola Norm. Sup. Pisa}, {\bf 22}, page 607--694, 1968.

\bibitem{bass}Bass, R.: ``Diffusions and Elliptic Operators'', Springer Verlag,  1998.


\bibitem{bolt-szn-1}Bolthausen, E., Sznitman, A.-S.: ``On the static and dynamic points of view for certain random walks in random environment'', {\it Methods and Applications of Analysis}, {\bf 9}(3), page 345--376, 2002.

\bibitem{bolt-szn}Bolthausen, E., Sznitman, A.-S.: ``Ten Lectures on Random Media'', DMV-Lectures, volume 32, Birkh\"auser, Basel, 2002. 

\bibitem{bolt-szn-zeit}Bolthausen, E., Sznitman, A.-S., Zeitouni, O.: ``Cut points and diffusive random walks in random environment'', {\it Ann. I. H. Poincar\'e}, PR 39(3), page 527--555, 2003.

\bibitem{com-zeit}Comets, F., Zeitouni, O.: ``A law of large numbers for random walks in random mixing environments'', {\it Ann. Probab.}, {\bf 32}(1B), page 880--914, 2004. 

\bibitem{deMasi}De Masi, A., Ferrari, P.A., Goldstein, S., Wick, W.D.: ``An invariance principle for reversible Markov processes. Applications to random motions in random environments'', {\it J. Statist. Phys.}, {\bf 55}, page 787--855, 1989.

\bibitem{gil-tru}Gilbarg D., Trudinger N.S.: ``Elliptic Partial Differential Equations of the Second Order'', Springer Verlag, 1998.

\bibitem{illin}Il'in, A.M., Kalashnikov, A.S., Oleinik, O.A.: ``Linear equations of the second order of parabolic type'', {\it Russian Math. Surveys}, {\bf 17}(3), page 1--143, 1962.

\bibitem{kalikow}Kalikow, S.A.: ``Generalized random walk in a random environment'', {\it Ann. Probab.}, {\bf 9}, page 753--768, 1981.

\bibitem{kar-shr}Karatzas, I., Shreve, S.: ``Brownian Motion and Stochastic Calculus'', Second Edition, Springer Verlag, 1991.

\bibitem{kip-Var}Kipnis, C., Varadhan, S.R.S.: ``A central limit theorem for additive functionals of reversible Markov processes and applications to simple exclusions'', {\it Commun. Math. Phys.}, {\bf 104}, 1--19, 1986.

\bibitem{kom}Komorowski, T.: ``Stationarity of Lagrangian velocity in compressible environments'', {\it Comm. Math. Phys.}, {\bf 228}(3), page 417--434, 2002.

\bibitem{kom-krupa02}Komorowski, T., Krupa, G.: ``On the existence of invariant measure for Lagrangian velocity in compressible environments'', {\it J. Statist. Phys.}, {\bf 106}(3-4), page 635--651, 2002.

\bibitem{kom-krupa}Komorowski, T., Krupa, G.: ``On stationarity of Lagrangian observations of passive tracer velocity in a compressible environment'', {\it Ann. Appl. Prob.}, to appear in Nov. 2004.

\bibitem{kom-olla-01}Komorowski, T., Olla, S.: ``On homogenization of time-dependent random flows'',  {\it Probab. Theory relat. Fields}, {\bf 121}(1), page 98--116, 2001.

\bibitem{kom-olla-03}Komorowski, T., Olla, S.: ``Invariant measures for passive tracer dynamics in Ornstein-Uhlenbeck flows'', {\it Stoch. Proc. Appl.}, {\bf 105}, page 139--173, 2003.

\bibitem{kozlov}Kozlov, S.M.: ``The method of averaging and walks in inhomogeneous environments'', {\it Russian Math. Surveys}, {\bf 40},  page 73--145, 1985.

\bibitem{landim-olla-yau}Landim, C., Olla, S., Yau, H.T.: ``Convection-diffusion equation with space-time ergodic random flow'', {\it Probab. Theory relat. Fields}, {\bf 112}, page 203--220, 1998.

\bibitem{lyons} Lyons, T.J., Zheng, W.A.: ``On conditional diffusion processes'', Proc. Roy. soc. Edinburgh Sect.A, {\bf 115}(3--4), page 243--255, 1990.

\bibitem{molchanov}Molchanov, S.A.: ``Lectures on random media'', {\it Lecture Notes in Math.}, volume 1581, page 242--411, Springer Verlag, 1994.

\bibitem{oel}Oelschl\"ager, K.: ``Homogenization of a diffusion process in a divergence-free random field'', {\it Ann. Probab.}, {\bf 16}(3), page 1084--1126, 1988.

\bibitem{olla94}Olla, S.: ``Homogenization of diffusion processes in random fields'', Ecole Doctorale, Ecole Polytechnique, Palaiseau, 1994.

\bibitem{olla01}Olla, S.: ``Central limit theorems for tagged particles and for diffusions in random environment''. In: ``Milieux Al\'eatoires'', Panoramas et Synth\`eses, Num\'ero 12, Soci\'et\'e Math\'ematique de France, 2001.

\bibitem{papa}Papanicolaou, G., Varadhan, S.R.S.: ``Diffusion with random coefficients'', {\it Statistics and probability: essays in honor of C.R. Rao}, G. Kallianpur, P.R. Krishnajah, J.K. Gosh, eds., North Holland, Amsterdam, page 547--552, 1982.

\bibitem{ras}Rassoul-Agha, F.: ``The point of view of the particle on the law of large numbers for random walks in a mixing random environment'', {\it Ann. Probab.}, {\bf 31}(3), page 1441--1463, 2003.

\bibitem{rev-yor}Revuz, D., Yor, M.: ``Continuous Martingales and Brownian Motion'', 3rd Edition, Springer Verlag, Berlin, 1999.

\bibitem{shen}Shen, L.: ``On ballistic diffusions in random environment'', {\it Ann. I. H. Poincar\'e}, PR 39(5), page 839--876, 2003.

\bibitem{shen-add}Shen, L.: Addendum to  ``On ballistic diffusions in random environment'', {\it Ann. I. H. Poincar\'e}, PR 40(3), page 385--386, 2004.

\bibitem{stroock}Stroock, D.: ``Diffusion semigroups corresponding to uniformly elliptic divergence form operators'', {\it Lecture Notes in Math.}, volume 1321, page 316--347, Springer Verlag, Berlin, 1988.

\bibitem{szn00}Sznitman, A.-S.: ``Slowdown estimates and central limit theorem for random walks in random environment'', {\it J. Eur. Math. Soc.}, {\bf 2},  page 93--143, 2000.

\bibitem{szn01}Sznitman, A.-S.: ``On a class of transient random walks in random environment'',  {\it Ann. Probab.},  {\bf 29}(2), page 723--764, 2001. 

\bibitem{szn02}Sznitman, A.-S.: ``An effective criterion for ballistic behavior of random walks in random environment'', {\it Probab. Theory relat. Fields}, {\bf 122}(4), page 509--544, 2002.

\bibitem{szn03}Sznitman, A.-S.: ``On new examples of ballistic random walks in random environment'', {\it Ann. Probab.}, {\bf 31}(1), page 285--322, 2003.

\bibitem{szn04}Sznitman, A.-S.: ``Topics in random walks in random environment'', {\it ICTP Lecture Notes Series}, Volume XVII: School and Conference on Probability Theory, May 2004.

\bibitem{szn-zer}Sznitman, A.-S., Zerner, M.P.W.: ``A law of large numbers for random walks in random environment'', {\it Ann. Probab.}, {\bf 27}(4), page 1851--1869, 1999.

\bibitem{zeit}Zeitouni, O.: ``Random Walks in Random Environment'', {\it Lecture Notes in Mathematics}, volume 1837, page 190--312, Springer, 2004.

\end{thebibliography}
\end{document}